\documentclass[twocolumn]{autart}
\usepackage{amsmath,amssymb,pstricks,color,dsfont}
\usepackage{graphicx}
\usepackage{natbib}
\usepackage{psfrag,graphicx,epsfig}
\usepackage{algorithm,algorithmic,xspace}
\usepackage{multirow}
\allowdisplaybreaks

\newcommand{\bremark}{\begin{remark} \begin{rm} }
\newcommand{\eremark}{ \end{rm} \rule{1mm}{2mm}
\end{remark} }
\newcommand{\btheorem}{\begin{theorem} \begin{rm} }
\newcommand{\etheorem}{ \end{rm} \rule{1mm}{2mm}
\end{theorem} }
\newcommand{\blemma}{\begin{lemma} \begin{rm} }
\newcommand{\elemma}{ \end{rm} \rule{1mm}{2mm}
\end{lemma} }
\newcommand{\bcorollary}{\begin{corollary} \begin{rm} }
\newcommand{\ecorollary}{ \end{rm} \rule{1mm}{2mm}
\end{corollary} }
\newcommand{\bdefinition}{\begin{definition}\begin{rm} }
\newcommand{\edefinition}{ \end{rm} \rule{1mm}{2mm}
\end{definition} }
\newcommand{\bproposition}{\begin{proposition} \begin{rm} }
\newcommand{\eproposition}{ \end{rm} \rule{1mm}{2mm}
\end{proposition} }
\newcommand{\bexample}{\begin{example} \begin{rm} }
\newcommand{\eexample}{ \end{rm} \rule{1mm}{2mm}
\end{example} }
\newcommand{\basm}{\begin{assumption} \begin{rm}}
\newcommand{\easm}{\end{rm} %\rule{1mm}{2mm}
\end{assumption}}

\newtheorem{theorem}{\bf Theorem}[section]
\newtheorem{lemma}{\bf Lemma}[section]
\newtheorem{definition}{\bf Definition}[section]
\newtheorem{remark}{\bf Remark}[section]
\newtheorem{corollary}{\bf Corollary}[section]
\newtheorem{proposition}{\bf Proposition}[section]
\newtheorem{example}{\bf Example}[section]
\newtheorem{assumption}{\bf Assumption}[section]

\newcommand\oprocendsymbol{\hbox{$\bullet$}}
\newcommand\oprocend{\relax\ifmmode\else\unskip\hfill\fi\oprocendsymbol}

\let\oldnl\nl% Store \nl in \oldnl
\newcommand{\nonl}{\renewcommand{\nl}{\let\nl\oldnl}}% Remove line number for one line

\date{}

\begin{document}

\begin{frontmatter}

\title{On privacy preserving data release of linear dynamic networks\thanksref{footnoteinfo}}

\thanks[footnoteinfo]{This work was partially supported by the grants ARO W911NF-13-1-0421 (MURI), NSF CNS-1505664, and NSF CAREER ECCS-1846706.}

\author[PSU]{Yang Lu}\ead{yml5046@psu.edu},
\author[PSU]{Minghui Zhu}\ead{muz16@psu.edu}

\address[PSU]{School of Electrical Engineering and Computer Science, Pennsylvania State University, 201 Old Main, University Park, PA, 16802 USA}


\begin{abstract}
Distributed data sharing in dynamic networks is ubiquitous. It raises the concern that the private information of dynamic networks could be leaked when data receivers are malicious or communication channels are insecure. In this paper, we propose to intentionally perturb the inputs and outputs of a linear dynamic system to protect the privacy of target initial states and inputs from released outputs. We formulate the problem of perturbation design as an optimization problem which minimizes the cost caused by the added perturbations while maintaining system controllability and ensuring the privacy. We analyze the computational complexity of the formulated optimization problem. To minimize the $\ell_0$ and $\ell_2$ norms of the added perturbations, we derive their convex relaxations which can be efficiently solved. The efficacy of the proposed techniques is verified by a case study on a heating, ventilation, and air conditioning system.
\end{abstract}

\end{frontmatter}

\section{Introduction\label{section introduction}}

Recently, information and communications technologies are increasingly integrated with control systems in the physical world. It has been stimulating the rapid emergence of cyber-physical systems (CPS). CPS consists of a large number of geographically dispersed entities and thus distributed data sharing is necessary to achieve network-wide goals. However, distributed data sharing also raises the significant concern that the private or confidential information of legitimate entities could be leaked to malicious entities. Privacy has become an issue of high priority to address before certain CPS can be widely deployed. For example, the current absence of accepted solutions to tackle privacy concerns caused a deadlock in the mandatory deployment of smart meters in the Netherlands because of the common belief that smart metering is necessarily privacy-invasive \citep{Cavoukian:tech12}. In 2010, California's new law on smart meter privacy indicated strong demands to protect the privacy of end-users' energy consumption data \citep{CA-Bill:10}.

%\begin{table*}[!t]
%\scriptsize
%\caption{Summary of our results}
%%\resizebox{\textwidth}{!}{
%\centering
%\begin{tabular}{|c|c|c|c|}
%  \hline & \text{non-convexity}
%  & \text{computational hardness} & \text{approach}
%\tabularnewline
%  \hline
%  $\mathbb{P}_0 \;(\|\cdot\|_0)$ & \text{objective function + constraint set} & \text{provably NP-hard} & $\ell_1$ + \text{nuclear\;\;norm} + \text{PSD}
%\tabularnewline
%  \hline
%  $\mathbb{P}_2 \;(\|\cdot\|_2)$ & \text{constraint set} & \text{N/A} & \text{SVD}
%\tabularnewline
%  \hline
%\end{tabular}
%%}
%\label{t2-table1}
%\end{table*}
%\normalsize

In the security community, several notions have been used to define privacy.
%, e.g., differential privacy \cite{CD:06,CD-AR:2014}, mutual information \cite{LS-SRR-HVP:2013,CES:1948}, semantic security \cite{OG:2004,SG-SM:1984,JK-YL:2007}, perfect secrecy \cite{Shamir,Shannon}, $k$-anonymity \cite{PS-LS:1998,LS:2002}, $\ell$-diversity \cite{CCA-PSY:2008,AM-DK-JG-MV:2007,EZ-LG:2011} and $t$-closeness \cite{NL-TL-SV:2007}.
In particular, differentially private schemes add random noises into each individual's data in such a way that, with high probability, the participation of the individual cannot be inferred by an adversary, who can access arbitrary auxiliary information, via released data \citep{CD-AR:2014}.
%Differential privacy protects individuals from additional harm caused by data release.
Mutual information \citep{LS-SRR-HVP:2013} requires explicit statistical models of source data and auxiliary/side information and quantifies average uncertainties about source data conditioned on revealed data.
%In order to hide source data, data owners aim to minimize the mutual information between source data and revealed data. In contrast to differential privacy, mutual information requires explicit statistical models of source data and auxiliary/side information. Another widely used privacy notion is semantic security \cite{OG:2004,SG-SM:1984,JK-YL:2007}.
Semantic security \citep{SG-SM:1984} requires that no additional information about a plaintext can be inferred using its ciphertext by any probabilistic polynomial-time algorithm. Perfect secrecy \citep{Shannon} is stronger than semantic security in that it assumes that the adversary has unlimited computing power. Additionally, $k$-anonymity \citep{PS-LS:1998} protects identity privacy by requiring that each group of records that share the same values for the quasi-identifiers (e.g., age, gender, zip code) must include at least $k$ records. The notion of $\ell$-diversity \citep{AM-DK-JG-MV:2007} extends $k$-anonymity to protect attribute privacy by requiring that there is adequate diversity in each sensitive attribute. The notion of $t$-closeness \citep{NL-TL-SV:2007} further refines $\ell$-diversity by taking into account side information of \emph{a priori} distributions of the attributes.

%In \cite{SH-UT-GJP:2016}, the authors pointed out that differential privacy is event based.

In the control and CPS communities, differential privacy has been adopted to Kalman filtering \citep{JLN-GJP:2014}, consensus \citep{ZH-SM:2012} and optimization \citep{EN-PT-JC:2016ACC,HZ-YS-PC-JC:2016}; mutual information has been used as a privacy metric in the applications of smart grid \citep{SH-UT-GJP:2016} and stochastic control systems \citep{PV-JY-PP:2015}; semantic security and perfect secrecy have been employed in secure multiparty computation \citep{YL-MZ:NECSYS2015Privacy} and homomorphic encryption \citep{YL-MZ:2018Automatica}.

\emph{Contributions.} This paper considers a linear dynamic network where a set of agents are physically coupled. An external data requester requires the agents to release system outputs in real time. The agents aim to prevent the data requester from inferring initial states and past inputs through the released outputs. We define that the privacy is protected if the data requester has infinite uncertainty on each of its target entries after observing the released outputs. Our uncertainty-based privacy notion extends $\ell$-diversity from the discrete-valued setting to the continuous-valued setting. Please refer to Section \ref{privacy notion section} and the appendix for the justification of our privacy definition.
%We realize the privacy requirement by a rank deficiency constraint and
We propose a protection scheme where the agents intentionally perturb the inputs and outputs such that (i) privacy is protected;
%(i) the rank deficiency constraint is satisfied;
(ii) system controllability is maintained; and (iii) the cost induced by the perturbations is minimized. We investigate two cases of the cost function:\\
(1) The sparsity of the perturbations is maximized, i.e., the $\ell_0$ norm of the added perturbations is minimized.\\
(2) The utility of the released outputs is maximized, i.e., the $\ell_2$ norm of the added perturbations is minimized.

We first analyze the computational complexity of the formulated optimization problem.
%induced by the non-convexity of the controllability and the rank deficiency constraints and the $\ell_0$ norm in the objective function.
We then derive a semidefinite program relaxation for the $\ell_0$ norm minimization by adopting the $\ell_1$ norm heuristic, the nuclear norm heuristic and a positive semidefinite condition. For the $\ell_2$ norm minimization, by using the tool of singular value decomposition, we provide a computationally more efficient method which can return an analytic feasible solution.
%The obtained results are summarized in Table \ref{t2-table1}.
Finally, the efficacy of the developed techniques is verified by a case study on a heating, ventilation, and air conditioning (HVAC) system.

%Similar to $\ell$-diversity \cite{NL-TL-SV:2007}, our privacy notation is vulnerable to skewness attacks. An interesting future work is to extend $t$-closeness, a refinement of $\ell$-diversity against skewness attacks, from the discrete-valued settings to our continuous-valued settings. Please refer to Section \ref{privacy notion section} for a detailed discussion.

%A preliminary version of this paper is presented in \cite{MZ-YL:ACC2015}. In comparison, in the current paper, the privacy definition is refined; a thorough analysis of computational intractability is provided; the relaxation for the $\ell_0$ minimization problem is studied; the $\ell_2$ minimization is extended for general rank deficiency constraints; a case study on an HVAC system is provided. Due to space limitation, the detailed introduction to $\ell$-diversity is omitted in the current paper and included in the appendix of the complete version \cite{YL-MZ-datareleaseAutomatica18-FV}.

A preliminary version of this paper is presented in \citep{MZ-YL:ACC2015}. Compared with \citep{MZ-YL:ACC2015}, in the current paper, the privacy definition is refined; a thorough analysis of computational complexity is provided; the relaxation for the $\ell_0$ minimization problem is developed; the $\ell_2$ minimization is extended to allow for general rank deficiency constraints; a case study on an HVAC system is provided.

\textbf{Notations and notions}. For any $k\in\mathbb{N}$, denote by $x_{[0,k]}$ the sequence of $\{x_0,\cdots,x_k\}$. The induced $\ell_1$ norm, $\ell_2$ norm and Frobenius norm of matrix $M$ are denoted by $\|M\|_1$, $\|M\|_2$ and $\|M\|_F$, respectively. $\|M\|_0$ denotes the number of nonzero entries of matrix $M$ and is referred to as its $\ell_0$ norm. $\|M\|_*$ denotes the sum of the singular values of matrix $M$ and is referred to as its nuclear norm.
%It is well-known that, for a matrix $M$ and a vector $x$ with compatible dimensions, $\|M x\|_2 \leq \|M\|_2\|x\|_2$, and for two matrices $M_1$ and $M_2$ with compatible dimensions, $\|M_1 M_2\|_2 \leq \|M_1\|_2\|M_2\|_2$.
For a column vector $w\in\mathbb{R}^n$, define a quantity $|\cdot|_{\min}$ as $|w|_{\min}=\min_{\ell\in\{1,\cdots,n\}}|w_\ell|$.
%Then, for any $w,w'\in\mathbb{R}^n$, $\|w-w'\|_{\min}$ is the entry-wise smallest distance between $w$ and $w'$.
$\mathbb{S}^{n}$ denotes the set of real symmetric matrices of size $n$. Denote by $M^\dag$ and $M^{-1}$ the pseudo inverse and inverse of matrix $M$, respectively.
%If $M$ is square and invertible, then the pseudo inverse reduces to the inverse, denoted by $M^{-1}$.
The notation $M=[M_{ij}]$ means that $M$ is a matrix for which the entry at the position of the $i$-th row and the $j$-th column is $M_{ij}$, or the $ij$-th block is $M_{ij}$ if $M_{ij}$ itself is a matrix. For a vector $x$, $x=[x_i]$ means that the $i$-th entry of $x$ is $x_i$.
%The notation ${\rm diag}\{M\}$ denotes the block diagonal matrix for which the sub-matrices of the diagonal blocks are $M$ and those of the off-diagonal blocks are zero matrices.
${\rm Tr}(M)$, ${\rm rank}(M)$ and ${\rm det}(M)$ denote the trace, the rank and the determinant of matrix $M$, respectively. 
%Denote by $\circ$ the operator of the Hadamard product of two matrices of the same size, i.e., given matrices $A=[A_{ij}]\in\mathbb{R}^{m\times n}$ and $B=[B_{ij}]\in\mathbb{R}^{m\times n}$, $T=[T_{ij}]=A\circ B$ means $T\in\mathbb{R}^{m\times n}$ and $T_{ij}=A_{ij}B_{ij}$. $1_{[{\rm condition}]}$ denotes the indicator such that $1_{[{\rm condition}]}=1$ if ${\rm condition}$ is true and $1_{[{\rm condition}]}=0$ if ${\rm condition}$ is false. 
Given a matrix $M$, denote by ${\rm vec}(M)$ the column vector consisting of the entries of $M$.
%For all positive integer $\ell$, $[\ell]$ is the power set of $\{1,\cdots,\ell\}$.
$\textbf{0}_n$ is the column vector with $n$ zeros. $\textbf{0}_{m\times n}$ denotes the $m\times n$ matrix where all entries are zeros. $I_n$ denotes the identity matrix of size $n$. In computational complexity, a polynomial-time algorithm is an algorithm performing its task in a number of steps bounded by a polynomial expression in the size of the problem input. NP (nondeterministic polynomial time), one of the most fundamental complexity classes, is the set of decision problems where the ``yes''-instances can be accepted in polynomial time by a nondeterministic Turing machine. A decision problem $\mathcal{D}$ is NP-hard if for every problem $\mathcal{D}'$ in NP, there is a polynomial-time reduction from $\mathcal{D}'$ to $\mathcal{D}$ \citep{JVL:1990}.

\section{Problem Statement} \label{sec:Problem}

This section introduces the system model, the adversary model and the privacy notion adopted in this paper.

\subsection{Network model}

%\subsubsection{Physical dynamics}

Consider an interconnected dynamic network of $V = \{1,\cdots,N\}$ where the physical dynamics of agent~$i$ are described by the following linear discrete-time system:
\begin{align}
x_i(k+1) &= \bar A_{ii} x_i(k) + \sum\nolimits_{j\in \mathcal{N}_i}\bar A_{ij}x_j(k) + \bar B_i u_i(k)\nonumber\\
y_i'(k) &= \bar G_i' x_i(k) + \bar H_i' u_i(k).
\label{e2}
\end{align}
In \eqref{e2}, $x_i(k)\in \mathbb{R}^{n_i}$, $u_i(k)\in \mathbb{R}^{p_i}$ and $y_i'(k)\in\mathbb{R}^{l_i}$ are the state, input and output of agent $i$ at time instant $k$, respectively,
%$\bar A_{ij}\in\mathbb{R}^{n_i\times n_j}$, $\bar B_i\in\mathbb{R}^{n_i\times p_i}$, $\bar G_i'\in\mathbb{R}^{l_i\times n_i}$, $\bar H_i'\in\mathbb{R}^{l_i\times p_i}$,
and $\mathcal{N}_i\subseteq V\setminus\{i\}$ is the set of agents whose states affect the state of agent $i$.
%Denote by $\GG = (V,E)$ the digraph representing the physical interconnections of the agents, where $E$ is the set of edges such that $(i,j)\in E$ if $j\in\NN_i$.
The collection of states and outputs in \eqref{e2} can be compactly written as follows:
\begin{align}
\label{state eq}
x(k+1) &= \bar A x(k) + \bar B u(k)\\
\label{output eq}
y'(k) &=\bar G'x(k)+\bar H'u(k)
\end{align}
where $x(k) = [x_i(k)]\in\mathbb{R}^n$, $u(k)=[u_i(k)]\in\mathbb{R}^p$ and $y'(k)=[y_i'(k)]\in\mathbb{R}^l$,
%$\bar A\in\mathbb{R}^{n\times n}$, $\bar B\in\mathbb{R}^{n\times p}$, $\bar G'\in\mathbb{R}^{l\times n}$ and $\bar H'\in\mathbb{R}^{l\times p}$,
with $n=\sum_{i\in V}n_i$, $p=\sum_{i\in V}p_i$ and $l=\sum_{i\in V}l_i$. The matrices $\bar A$, $\bar B$, $\bar G'$ and $\bar H'$ are system parameters known to the agents. In this paper, we assume that each agent is allowed to communicate with any other agent and measure all the entries of its own state. The activation of a communication link or a sensor induces certain cost.
%Assume that $n\geq q\geq p$.
%\subsubsection{Communication and sensing}
%Each agent can communicate with some other agents. The inter-agent communication topology is denoted by digraph $\mathcal{G}^C=(V,E^C)$ where $(i,j)\in E^C$ if agent $j$ can send messages to agent $i$. Let $S^x$ (resp. $S^u$) be the set of $(i,\ell)$
%with $i\in V$ and $\ell\in\{1,\cdots,n_i\}$ (resp. $\ell\in\{1,\cdots,p_i\}$)
%such that $x_{i\ell}$ (resp. $u_{i\ell}$) can be measured by agent $i$. Hence, $S^x$ and $S^u$ represents the set of sensor locations and $\bar G'$ and $\bar H'$ in \eqref{output eq} are restricted by $S^x$ and $S^u$, respectively. {\color{blue}Specifically, if $(i,\ell)\notin S^x$ (resp. $(i,\ell)\notin S^u$), then the $\ell$-th column of $\bar G_i'$ (resp. $\bar H_i'$) is a zero column.}

\subsection{Adversary model}

There is an external data requester who requests a set of linear combinations of the agents' individual outputs. Specifically, the data requester determines a constant matrix $\Pi\in\mathbb{R}^{q\times l}$ and tells its valuation to a data aggregator. Each agent $i$ measures its output $y_i'(k)$ and sends it to the data aggregator, who then computes $y(k)=\Pi y'(k)$ and sends $y(k)$ to the data requester. Hence, the data received by the data requester is:
\begin{align}
y(k)= \Pi y'(k)=\bar Gx(k)+\bar Hu(k)
\label{e9}
\end{align}
where $\bar G=\Pi \bar G'\in\mathbb{R}^{q\times n}$ and $\bar H=\Pi \bar H'\in\mathbb{R}^{q\times p}$.  
%An external data requester requires each agent $i$ to send him a weighted output $\Pi_iy_i'(k)$, where $\Pi_i\in\mathbb{R}^{q\times l_i}$ is the weight matrix and $y_i'(k)$ is agent $i$'s output in \eqref{output eq}. The data requester then sums the received data. $\sum\nolimits_{i=1}^N\Pi_iy_i'(k)=$ 
%{\color{blue}Specifically, each agent $i$ sends $\Pi_iy_i'(k)$ to the data requester, where $\Pi_i\in\mathbb{R}^{q\times l_i}$ is the sub-matrix of $\Pi$ corresponding to $y_i'$. The data requester then computes the sum $y(k)=\Pi y'(k)=\sum_{i=1}^N\Pi_iy_i'(k)$.} 
The agents are unaware of how the released data will be used. In the rest of the paper, we use $(A,B,G,H)$ to represent arbitrary system matrices, while $(\bar A,\bar B,\bar G,\bar H)$ specifically represent system \eqref{state eq} and \eqref{e9}.

%\begin{figure}[H]
%begin{center}
%\includegraphics[width=1\linewidth]{figures/agent_response}
%\caption{Data release of linear dynamic networks}
%\label{agent_response}
%\end{center}
%\end{figure}

The data requester is assumed to be semi-honest
%(or passive or honest-but-curious)
\footnote{Semi-honest adversaries correctly follow the algorithm but attempt to use the received messages to infer private/confidential information of legitimate entities.} \citep{YL-BP:09}, and aims to exploit $y(k)$ to infer some entries of $x(0)$ and $\{u(k)\}$. Our problem model is motivated by several practical scenarios, e.g., smart building and load monitoring in smart grid. For both of these two applications, the physical dynamics can be approximated by a linear dynamic model \citep{CSK-JIP-MP-JB:2014,IM-YT:2014}. In smart building, a system operator (data requester) uses temperature data of individual rooms (agents) to monitor working comfortability and energy usage conditions. At the same time, it may aim to infer occupancy data from room temperature information, and by the occupancy data, it might be able to derive the location traces of individual occupants \citep{ML-DM-SW:2010}. In smart grid, a utility company (data requester) collects power consumption data stored at local smart meters of power consumers (agents) to monitor power usage conditions. Meanwhile, it may target to infer power load profiles of individual consumers from aggregated home power consumption information \citep{SM-PM-WA:09}.%\footnote{It is equivalent to assume that the data requester is benign but adversaries can eavesdrop the messages sent to the data requester.}
We assume that the data requester is aware of the matrices $\bar A$, $\bar B$, $\bar G'$ and $\bar H'$. This assumption models the auxiliary/side information of the adversary. In this paper, we assume that all the agents in $V$ and the data aggregator are benign, i.e., they will not use their observed information to infer other agents' private data.

%This assumption makes privacy protection challenging, because it considers the worst case where the adversary data requester knows everything about the system. In cryptography, this is called Kerckhoffs' principle or ``the enemy knows the system'' principle \cite{Shannon} and has been widely embraced by cryptographers.

%One of our future directions is to relax this assumption \cite{KA-WQ-YH:97}.

%Semi-honest adversary model has been broadly used in various applications, e.g., secure multiparty computation, linear programming, dataset process and consensus \cite{RC-ID-JBN:2015,JD-FK:2011,MJF-KN-BP,ZH-SM:2012}.
%In addition, it is reasonable to assume that the monitoring operator is a semi-honest inference attacker.\margin{inference attacks have been removed from the introduction} As discussed in the introduction, inference attacks aim to steal private information for subsequent direct attacks to dynamic networks. So inference attackers are not willing to trigger intrusion detection and get caught before successfully stealing private information.

%The above model can capture the following two scenarios:
%\begin{itemize}
%\item System~\eqref{e8} is a controlled system and $D(k)$ is the control input.
%\item System~\eqref{e8} is an autonomous system and $GD(k)$ represents the unknown
%\end{itemize}

\subsection{Privacy notion\label{privacy notion section}}

We next introduce the privacy notion adopted in this paper. The data requester aims to infer the values of partial (could be all) entries of the initial state $x(0)$ and the input sequence $\{u(k)\}$. We call those entries as target entries. The remaining entries of $x(0)$ and $\{u(k)\}$ are called nontarget entries. Denote by $x^t(0)$ and $x^n(0)$ (resp. $u^t$ and $u^n$) the column vectors of the target and nontarget entries of $x(0)$ (resp. $u$), respectively. Denote by $d_x^t$, $d_u^t$, $d_x^n$ and $d_u^n$ the dimensions of $x^t(0)$, $u^t$, $x^n(0)$ and $u^n$, respectively. It holds that $d_x^t+d_x^n=n$ and $d_u^t+d_u^n=p$. Denote by $x_\ell^t(0)$ (resp. $u_\ell^t$, $x_\ell^n(0)$ and $u_\ell^n$) the $\ell$-th entry of $x^t(0)$ (resp. $u^t$, $x^n(0)$ and $u^n$). For a target entry $u_{\ell}^t$, we consider that it is protected if and only if $u_{\ell}^t(k)$ is protected for any $k\in\mathbb{N}$. In other words, if the value of $u_{\ell}^t(k)$ for one time instant $k$ is disclosed to the data requester, then we consider that the privacy of $u_{\ell}^t$ is compromised.

Given system matrices $(A,B,G,H)$,
%by Section 4.2.2 of \cite{Chen},
for each time instant $k$, the output $y(k)$ can be expressed as a linear combination of the entries of $x(0)$ and $u_{[0,k]}$:
\begin{align}
\label{standard output relation}
\!\!\!y(k)=GA^kx(0)+\sum\limits_{m=0}^{k-1}GA^{k-1-m}Bu(m)+Hu(k).
\end{align}
%Given system matrices $(A,B,G,H)$ and time instant $\kappa\in\mathbb{N}$, denote by $Y_{A,B,G,H,\kappa}$ the set of all possible output sequence up to time instant $\kappa$ generated by system $(A,B,G,H)$, i.e., $Y_{A,B,G,H,\kappa}=\{y(0),\cdots,y(\kappa):\exists x(0),u(0),\cdots,u(\kappa)\;{\rm s.t.}\;{\rm Eq.}\;\eqref{standard output relation}\;{\rm holds}\;\forall k=0,\cdots,\kappa\}$.
%By a rearrangement of the system matrices in \eqref{standard output relation}, we can express $y(k)$ as $y(k)=f_{A,B,G,H,k}^t(x^t(0),u^t(0),\cdots,u^t(k))+f_{A,B,G,H,k}^n(x^n(0),u^n(0),\cdots,u^n(k))$, where $f_{A,B,G,H,k}^t:\mathbb{R}^{d_x^t+(k+1)d_u^t}\to\mathbb{R}^{q}$ (resp. $f_{A,B,G,H,k}^n:\mathbb{R}^{d_x^n+(k+1)d_u^n}\to\mathbb{R}^{q}$) represents the linear function that maps $(x^t(0),u^t(0),\cdots,u^t(k))$ (resp. $(x^n(0),u^n(0),\cdots,u^n(k))$) to $y(k)$. Notice that $f_{A,B,G,H,k}^t$ and $f_{A,B,G,H,k}^n$ depends on $(A,B,G,H)$ and $k$.
Given system matrices $(A,B,G,H)$ and time instant $\kappa\in\mathbb{N}$, for any feasible output sequence $y_{[0,\kappa]}$, we define a set $\Delta_{A,B,G,H}(y_{[0,\kappa]})$ as:
\begin{align*}
&\Delta_{A,B,G,H}(y_{[0,\kappa]})=\{x^t(0),u^t_{[0,\kappa]}:\exists x^n(0),u^n_{[0,\kappa]},\\
&{\rm s.t.}\;y(k)={\rm right \text{-}hand \text{-}side\;of}\;\eqref{standard output relation},\forall k=0,\cdots,\kappa,\\
&{\rm with}\; x(0)\;{\rm the\;composition\;of\;}x^t(0)\;{\rm and}\;x^n(0)\;{\rm and}\;u(k)\\
&{\rm the\;composition\;of}\;u^t(k)\;{\rm and}\;u^n(k),\forall k=0,\cdots,\kappa\}.
\end{align*}
%\begin{align*}
%&\Delta_{A,B,G,H,\kappa}(y(0),\cdots,y(\kappa))\\
%&=\{x^t(0),u^t(0),\cdots,u^t(\kappa):\exists x^n(0),u^n(0),\cdots,u^n(\kappa),\\
%&\;\;\quad{\rm s.t.}\;y(k)=f_{A,B,G,H,k}^t(x^t(0),u^t(0),\cdots,u^t(k))\\
%&\;\quad+f_{A,B,G,H,k}^n(x^n(0),u^n(0),\cdots,u^n(k)),\;\forall k=0,\cdots,\kappa\},
%\end{align*}
The set $\Delta_{A,B,G,H}(y_{[0,\kappa]})$ includes all possible valuations of $\{x^t(0),u^t_{[0,\kappa]}\}$ that can generate $y_{[0,\kappa]}$ in \eqref{standard output relation}.
%Intuitively, given a sequence of observations $\{y(0),\cdots,y(\kappa)\}$, the larger the size of $\Delta_{A,B,G,H,\kappa}(y(0),\cdots,y(\kappa))$, the more uncertain the data requester is about the target entries. To define privacy,
%We next quantify the diameter of $\Delta_{A,B,G,H,\kappa}(y(0),\cdots,y(\kappa))$. Recall that, for any $w,w'\in\mathbb{R}^r$, $\|w-w'\|_{\min}$ is the entry-wise smallest distance between $w$ and $w'$; please refer to Section \ref{section notations}.
%For convenience of notation, for a column vector $w$ of proper dimension, we say $w\in\Delta_{A,B,G,H,\kappa}(y(0),\cdots,y(\kappa))$ if $w=[x^t(0)^T,u^t(0)^T,\cdots,u^t(\kappa)^T]^T$ and $\{x^t(0),u^t(0),\cdots,$ $u^t(\kappa)\}\in\Delta_{A,B,G,H,\kappa}(y(0),\cdots,y(\kappa))$.
We define the diameter of $\Delta_{A,B,G,H}(y_{[0,\kappa]})$ as:
\begin{align*}
&{\rm Diam}_{A,B,G,H}(y_{[0,\kappa]})=\sup\limits_{w,w'\in\Delta_{A,B,G,H}(y_{[0,\kappa]})}|w-w'|_{\min}.
\end{align*}
\begin{definition}%[Component-wise privacy]
\label{def: agent-wise privacy}
Given system matrices $(A,B,G,H)$, the privacy of $x^t(0)$ and $u^t$ is said to be protected if, for any $\kappa\in\mathbb{N}$, %$\Delta_{A,B,G,H,\kappa}(y(0),\cdots,y(\kappa))$ is path-connected and
${\rm Diam}_{A,B,G,H}(y_{[0,\kappa]})=\infty$
%for any $\{y(0),\cdots,y(\kappa)\}\in Y_{A,B,G,H,\kappa}$.
for any feasible output sequence $y_{[0,\kappa]}$.
\hfill\oprocend
\end{definition}

%Definition \ref{def: agent-wise privacy} states that the privacy of the target entries is protected if, given the data requester's observations, the size of the uncertainty on each target entry is infinite.
We next justify Definition \ref{def: agent-wise privacy} through comparisons with several popular existing notions in our problem setting.
%Please refer to Section \ref{justification} for the justification of Definition \ref{def: agent-wise privacy}.

%\subsection{Formulation justification\label{justification}}

%In this subsection, we justify our privacy notion through comparisons with several popular existing notions.

$\bullet$ Why not semantic security or perfect secrecy? These notions require that ``nothing is learned'' by the adversary from outputs. However, as pointed out in Section 2.2 of \citep{CD-AR:2014}, such ``nothing is learned'' definition cannot be adopted to applications in which the outputs have to be used to realize certain utility by the data user who is adversary, because such a strong privacy requirement intrinsically inhibits any meaningful data utility. In our problem, the adversary and the data user is the same entity, i.e., the data requester, and it has to accomplish certain analysis using the outputs.

$\bullet$ Why not mutual information metric?
%This notion provides a metric to quantify the amount of input information that can be inferred from outputs.
The usage of mutual information metric requires explicit statistical models of source data and auxiliary/side information \citep{LS-SRR-HVP:2013}. This requirement might be restrictive or even unrealistic for our problem as the inputs of the system may not follow any probabilistic distribution.

$\bullet$ Why not differential privacy?
%Differential privacy can quantify how much privacy can be achieved in a mathematically convenient and rigorous manner.
To achieve differential privacy, noises are persistently added to the released data via following, e.g., Gaussian and Laplace distributions. For control systems, such open-loop and persistent noise injection mechanisms could potentially deteriorate system performance.

$\bullet$ Our uncertainty-based privacy notion Definition \ref{def: agent-wise privacy}. Definition \ref{def: agent-wise privacy} is extended from the notion of $\ell$-diversity.
%As mentioned, the notion of $\ell$-diversity \cite{AM-DK-JG-MV:2007} is an extension of the notion of $k$-anonymity \cite{PS-LS:1998}.
In particular, $\ell$-diversity has been widely used in both application-centric research \citep{PMVK-MK:2012,NL-KD:2013} and formal privacy analysis \citep{ML-TT-AZ-ES:2014,NL-TL-SV:2007} on attribute privacy of discrete-valued tabular datasets. Besides academic studies, $k$-anonymity and $\ell$-diversity have also been popular in real world applications. As mentioned in page 16 of \citep{BM-PK-AH:2017}, $k$-anonymity has become a standard privacy notion in the industry. For example, the password manager 1Password has applied $k$-anonymity to protect the privacy of the customers' passwords \citep{JB:2018}. Recently, Google released a data loss prevention application programming interface (API) which supported $k$-anonymity and $\ell$-diversity \citep{CH:2017}.

Informally speaking, possessing $\ell$-diversity means that there are at least $\ell$ different values for each sensitive attribute of the dataset in the released table. A larger diversity indicates a larger uncertainty and thus the notion of diversity can be viewed as a measure of uncertainty. To make an analogy to $\ell$-diversity, in our problem, we can view each entry of $x^t(0)$ and $u^t$ as a sensitive attribute and require adequate diversity/uncertainty on it. In $\ell$-diversity, the diversity of discrete-valued sensitive attributes is defined by the number of different valuations for the attributes. In contrast, the target entries $x^t(0)$ and $u^t$ in our problem are continuous-valued and uncountable, which requires a new measure to quantify the diversity/uncertainty. In this paper, given system matrices $(A,B,G,H)$, the diversity/uncertainty is measured by the diameter of the set $\Delta_{A,B,G,H}(y_{[0,\kappa]})$.
%A larger diameter indicates a larger diversity/uncertainty on the target entries and an infinite diameter achieves the largest possible diversity/uncertainty.
%This is the rationale of Definition \ref{def: agent-wise privacy}.
%Hence, for our problem setting, we say that privacy is preserved if the diameter of $\Delta_{\hat A,\hat B,\hat G,\hat H,\kappa}(y(0),\cdots,y(\kappa))$ is infinite for any feasible output sequence $\{y(0),\cdots,y(\kappa)\}$ and any $\kappa\in\mathbb{N}$.
%In this paper, we measure the diversity/uncertainty of each target entry by the diameter of its admissible set given the data requester's observations $\{y(k)\}$ and we consider the case of infinite diversity/uncertainty to be adequate, for which we say that privacy is preserved.
Hence, Definition \ref{def: agent-wise privacy} extends the notion of $\ell$-diversity from the discrete-valued setting to the continuous-valued setting. A detailed introduction to $\ell$-diversity and extension to Definition \ref{def: agent-wise privacy} in our problem setting is given in the appendix.

Definition \ref{def: agent-wise privacy} is closely relevant to non-strong observability \citep{MLJH:1983} in control theory. Specifically, a dynamic system is not strongly observable if at least one entry of the initial states and input sequence is unobservable, i.e., cannot be uniquely determined. However, non-strong observability does not necessarily imply that all the target entries are unobservable. Definition \ref{def: agent-wise privacy} extends non-strong observability by explicitly ensuring such property. Using the language of control theory, Definition \ref{def: agent-wise privacy} can be equivalently stated as follows: Given system matrices $(A,B,G,H)$, the privacy of the target entries is said to be protected if no target entry is in the strongly observable subspace of system $(A,B,G,H)$. This uncertainty/unobservability-based privacy definition has been widely adopted in the control community; please see, e.g., \citep{YM-RMM:2017} and \citep{SP-SK-SS-APA:2014}, in which the initial state of a system is private if it is not in the observable subspace.
%The work \citep{YM-RMM:2017} focused on the characterization of the private/unobservable subspace for a given system, but did not study the problem of how to protect an initial state if it lies in the observable subspace of the system. The work \citep{SP-SK-SS-APA:2014} studied minimization of the dimension of the observable subspace. In contrast to \citep{YM-RMM:2017}, the current paper aims to protect initially insecure/observable target entries via perturbation design. In contrast to \citep{SP-SK-SS-APA:2014}, the current paper studies optimal perturbation design such that the strongly observable subspace of the perturbed system does not include any target entry, while system controllability is maintained and the perturbation cost is minimized.}

In discrete event systems, the notion of opacity has been widely used to define system state privacy \citep{YCW-SL:2014,BR-RC-SIM:2016,YJ-YCW-SL:2018}. 
%\citep{JD-PD-HM:2010,AS-CNH:2012,YCW-SL:2014,XY-SL:2016,YJ-YCW-SL:2018}
%Recently, opacity has been applied to linear dynamic systems \citep{BR-RC-SIM:2016}.
%\citep{BR-RC-SIM:2016,BR-RC-SIM:2016Allerton,BR-RC-SIM:2017}.
A system is opaque if for every secret-induced behavior, there exists a non-secret-induced behavior that generates identical observations. The notion of opacity is similar to our privacy notion in spirit. However, the privacy objectives are different. In opacity-based works, the privacy objective is to ensure that the adversary cannot determine from the observations whether or not the system state belongs to a predefined secret set. 
%In other words, opacity-based works protect state membership privacy. 
In contrast, Definition \ref{def: agent-wise privacy} aims to protect data privacy such that for any valuation of any target entry, the adversary has infinite uncertainty from the observations on the value of the target entry.

\emph{Advantages.} Compared with semantic security and perfect secrecy, our notion is weaker than the ``nothing is learned'' requirement and allows for meaningful data utility. Compared with mutual information, our notion does not require any statistical model for system states, inputs and outputs. Compared with differential privacy, our notion does not require using persistent perturbations. Please refer to Remark \ref{remark tradeoff} for more discussions.% also adopts perturbations for privacy preservation. However, by \eqref{e36}, it can be seen that the perturbations are added in a closed-loop fashion and diminishing as the system is stabilized. Since we formulate problem \eqref{OP agent-wise} such that the perturbed system remains controllable, one can design a feedback controller by the perturbed system matrices $(\bar A+\bar BK_{SS},\bar B(I_p+K_{SI}))$ to achieve perfect stability where the perturbations vanish at the equilibrium.

\emph{Limitation.}
%Differential privacy is able to neutralize any linkage attacks attempted with all forms and sources of auxiliary information (Section 2.3.2 of \cite{CD-AR:2014}).
A limitation of our privacy notion is that it does not take into account the scenario where the data requester has auxiliary information of some \emph{a priori} skewed distribution of $x(0)$ and $\{u(k)\}$. Note that $\ell$-diversity is also vulnerable to skewness attacks \citep{NL-TL-SV:2007}. For this case, instead of requiring infinite uncertainty on the target data items, the privacy goal should be that the posterior uncertainties after seeing the observations should be as close as possible to the \emph{a priori} uncertainties determined by the \emph{a priori} skewed distribution. This privacy goal extends $t$-closeness \citep{NL-TL-SV:2007}, which is a refinement of $\ell$-diversity, from discrete-valued settings to continuous-valued settings. We leave the study of the refined privacy goal as a future work.

\subsection{Privacy preserving data release\label{competitive privacy preserving section}}

To protect privacy, we propose to perturb the inputs and outputs such that the data requester cannot infer the target entries in the sense of Definition \ref{def: agent-wise privacy}. However, the perturbations should maintain certain system utilities, e.g., controllability. Throughout this paper, we assume that the original system $(\bar A,\bar B)$ is controllable and aim to maintain controllability of the perturbed system. These partially conflicting sub-objectives define the problem of \emph{privacy preserving data release}. In the remainder of the paper, we introduce our solutions of this problem. First, in Section \ref{sec: IOP}, we introduce our perturbation mechanism and the optimization formulation to solve for the optimal perturbation. In particular, we formulate an $\ell_0$ optimization that studies economy-privacy tradeoff and an $\ell_2$ optimization that studies utility-privacy tradeoff. After that, we show that the formulated optimization problems are hard to solve. In Section \ref{P0 section}, we first provide a further computational complexity result for the $\ell_0$ optimization problem and then derive a convex relaxation for it. A convex relaxation for the $\ell_2$ optimization problem is derived in Section \ref{Section P2}.

\section{Intentional input-output perturbations\label{sec: IOP}}

In this section, we first introduce a class of optimization problems to formulate intentional input-output perturbations. After that, we analyze the computational complexity of the formulated optimization problem.
%{\color{blue}In the next two sections, we introduce relaxation methods that can efficiently derive perturbed systems.}
%verify the computational hardness of the concerned problems.

%\begin{figure}[H]
%\begin{center}
%\includegraphics[width=1\linewidth]{figures/overall_scheme}
%\caption{Overall architecture of the proposed privacy preserving strategy {\color{red}(I am thinking to use either Fig. 2 or 3)}}
%\label{overall_scheme}
%\end{center}
%\end{figure}

%\begin{figure}[H]
%\begin{center}
%\includegraphics[width=1\linewidth]{figures/overall_scheme3}
%\caption{Overall architecture of the privacy preserving strategy}
%\label{overall_scheme}
%\end{center}
%\end{figure}

\subsection{Optimization formulation\label{OP formulation}}

To protect privacy, we propose the approach of \emph{intentional input-output perturbations}. Each agent $i$ intentionally perturbs its own input $u_i(k)$ and output $y_i'(k)$ by adding signals $\mu_i^\mu(k)\in\mathbb{R}^{p_i}$ and $\mu_i^y(k)\in\mathbb{R}^{l_i}$, respectively. The perturbations $\mu_i^u(k)$ and $\mu_i^y(k)$ are linear combinations of system states and inputs and given by: \begin{align}
\mu_i^u(k) &= \sum\nolimits_{j\in V}K_{ij}^{SS}x_j(k)+\sum\nolimits_{j\in V}K_{ij}^{SI}u_j(k)\nonumber\\
\mu_i^y(k) &= \sum\nolimits_{j\in V}K_{ij}^{OS}x_j(k)+\sum\nolimits_{j\in V}K_{ij}^{OI}u_j(k).
\label{e36}
\end{align}
%where $\mu_i^\mu(k)\in\mathbb{R}^{p_i}$ and $\mu_i^y(k)\in\mathbb{R}^{l_i}$.
%$K_{ij}^{SS}\in\mathbb{R}^{p_i\times n_j}$, $K_{ij}^{SI}\in\mathbb{R}^{p_i\times p_j}$, $K_{ij}^{OS}\in\mathbb{R}^{l_i\times n_j}$ and $K_{ij}^{OI}\in\mathbb{R}^{l_i\times p_j}$.
The superscript $SI$ means a perturbation from an input to a state. Other superscripts are defined analogously, with $O$ denoting output. Substituting the perturbations $\mu^u(k)=[\mu_i^u(k)]$ and $\mu^y(k)=[\mu_i^y(k)]$ into \eqref{state eq} and \eqref{e9} renders the following perturbed system:
\begin{align}
\label{perturbed system eq}
x(k+1) &= \bar A x(k) + \bar B(u(k) + \mu^u(k)) \nonumber\\
&= \hat A x(k)+\hat Bu(k)\\
\label{e12}
y(k) &= \Pi (\bar G' x(k) + \bar H'(u(k) + \mu^u(k)) + \mu^y(k))\nonumber\\
&=\hat Gx(k)+\hat Hu(k)
\end{align}
where $\hat A=\bar A+\bar B K_{SS}$, $\hat B=\bar B(I_p + K_{SI})$, $\hat G=\bar G+ \bar H K_{SS} + \Pi K_{OS}$ and $\hat H=\bar H + \bar H K_{SI} + \Pi K_{OI}$, with $K_{SS} = [K_{ij}^{SS}]\in\mathbb{R}^{p\times n}$, $K_{SI} = [K_{ij}^{SI}]\in\mathbb{R}^{p\times p}$, $K_{OS} = [K_{ij}^{OS}]\in\mathbb{R}^{l\times n}$ and $K_{OI} = [K_{ij}^{OI}]\in\mathbb{R}^{l\times p}$. Let $K=\left[\begin{array}{cc}
                               K_{SS} & K_{SI}\\
                               K_{OS} & K_{OI}
                             \end{array}\right]\in\mathbb{R}^{(p+l)\times(n+p)}$. In the rest of the paper, we use $(\hat A,\hat B,\hat G,\hat H)$ to specifically represent the perturbed system \eqref{perturbed system eq} and \eqref{e12}. The perturbation matrix $K$ is subject to two constraints:\\
%\textbf{(i)} The restrictions on the sensing and communication capabilities specified by $E^C$, $S^x$ and $S^u$.\\
\textbf{(i)} The perturbed system $(\hat A,\hat B)$ remains controllable.\\
\textbf{(ii)} The data requester cannot infer the target entries in the sense of Definition \ref{def: agent-wise privacy} from the outputs \eqref{e12}.

%The first constraint is captured by $\mathbb{K}\subseteq\mathbb{R}^{(p+l)\times(n+p)}$, which specifies the zero-nonzero structure of $K$. For example, the specification could be that $K\in\mathbb{K}$ if and only if $K$ is in the form of $\left[ {\begin{array}{*{20}{c}}
 %* &0& 0 \\
 %* &0& *\end{array}} \right]$, where $*$ means that the corresponding entry of $K$ could be any real number, while $0$ means that the corresponding entry of $K$ must be $0$ due to that no sensor and/or communication link can be installed at this position. This is analogous to pattern matrices studied in structured system theory \citep{CTL:1974}.
By adding $\mu^u$ to $u$ according to \eqref{e36}, the perturbed input $\hat u$ is $\hat u=u+\mu^u=K_{SS}x+(I_p+K_{SI})u$ and this actually changes system matrices $(\bar A,\bar B)$ to $(\hat A,\hat B)$. The controllability of $(\bar{A},\bar{B})$ does not guarantee that of $(\hat{A},\hat{B})$.
 %Hence, even if the original system $(\bar A,\bar B)$ is controllable, it is not guaranteed that there exists $u$ such that $\hat u$ in the above form can stabilize the system. For example, if $K_{SS}=\textbf{0}_{p\times n}$ and $K_{SI}=-I_p$, then the state equation of the perturbed system \eqref{perturbed system eq} reduces to $x(k+1)=\bar Ax(k)$. If $\bar A$ is not stabilizable, then for any $u$, the perturbed system diverges and hence the perturbed system is actually uncontrollable. Therefore, we need to reanalyze and construct the controllability with respect to the system matrices $(\hat A,\hat B)$ of the perturbed system.
 Denote by ${\mathcal{C}}(K_{SS},K_{SI})$ the controllability matrix of the perturbed system, i.e.,\\
 ${\mathcal{C}}(K_{SS},K_{SI})=[\hat B,\hat A\hat B,\cdots,\hat A^{n-1}\hat B]\\
 =[\bar B(I_{p}+K_{SI}),(\bar A+\bar BK_{SS})\bar B(I_{p}+K_{SI}),\\
 \quad\;\cdots,(\bar A+\bar BK_{SS})^{n-1}\bar B(I_{p}+K_{SI})]$.\\
%\begin{align*}
%{\mathcal{C}}(K_{SS},K_{SI})=&\left[ {\begin{array}{*{20}{c}}
%B(I_{p}+K_{SI})&(A+BK_{SS})B(I_{p}+K_{SI})
%\end{array}} \right.\\
%&\left. {\begin{array}{*{20}{c}}
%\cdots&(A+BK_{SS})^{n-1}B(I_{p}+K_{SI})
%\end{array}} \right].
%\end{align*}
The perturbed system is controllable if and only if ${\rm det}({\mathcal{C}}(K_{SS},K_{SI})({\mathcal{C}}(K_{SS},K_{SI}))^T)$ $>0$. Meanwhile, the agents aim to minimize the cost induced by the perturbations. This is captured by minimizing an objective function $c(K)$ determined later. All the above objectives are encoded in the following optimization problem:
\begin{align}
\label{OP agent-wise}
&\min\nolimits_{K\in\mathbb{R}^{(p+l)\times(n+p)}}c(K)\nonumber\\
&{\rm s.t.}%\;\Delta_{\hat A,\hat B,\hat G,\hat H,\kappa}(y(0),\cdots,y(\kappa)){\rm\; is\; path\text{-}connected}\nonumber\\
%&\quad\;\; {\rm and}
\,{\rm Diam}_{\hat A,\hat B,\hat G,\hat H}(y_{[0,\kappa]})=\infty,\forall\kappa\in\mathbb{N}\;{\rm and}\;{\rm feasible}\; y_{[0,\kappa]},\nonumber\\
%&\quad\;\; {\rm of\; system\;} (\hat A,\hat B,\hat G,\hat H),\nonumber\\
%&\quad\;\; \forall\{y(0),\cdots,y(\kappa)\}\in Y_{\hat A,\hat B,\hat G,\hat H,\kappa},\;\forall\kappa\in\mathbb{N},\nonumber\\
&\quad\;\; {\rm det}({\mathcal{C}}(K_{SS},K_{SI})({\mathcal{C}}(K_{SS},K_{SI}))^T) >0.
%&\quad\;\; \left[\begin{array}{cc}
%                               A+B K_{SS} & B(I_{p}+K_{SI})
%                             \end{array}
%                           \right]\;\;{\rm is\;\; controllable}.
\end{align}
%\begin{align}
%&\min\limits_{K\in\mathbb{K}}\|K\|_{\rm norm}\nnum\\
%&{\rm s.t.}\;\left[\begin{array}{cc}
%                               A+B K_{SS} & B(I_{p}+K_{SI}) \\
%                               \Pi(G'+H' K_{SS} + K_{OS}) & \Pi(H'+H' K_{SI} + K_{OI}) \\
%                             \end{array}
%                           \right]\nnum\\
%&\quad\quad{\rm is\;\;not\;\;strongly \;\;observable},\nnum\\
%&\quad\quad \left[\begin{array}{cc}
%                               A+B K_{SS} & B(I_{p}+K_{SI})
%                             \end{array}
%                           \right]\;\;{\rm is\;\; controllable}.
%\label{e3}
%\end{align}
In this paper, we study the following two representative cases of the cost function.

%\begin{itemize}
Problem $\mathbb{P}_0$: economy--privacy tradeoff. The added perturbations require communication and sensing. If one entry of $K^{SS}_{ji}$ or $K^{OS}_{ji}$ is nonzero, then agent $i$ needs to measure the corresponding entry of $x_i$ and sends it to agent $j$. If one entry of of $K^{SI}_{ji}$ or $K^{OI}_{ji}$ is nonzero, then agent $i$ needs to share its control $u_i$ with agent $j$. Recall that activation of communication links and sensors induces some cost. 
    %If one element of $K^{SS}_{ij}$ or $K^{SO}_{ij}$ (resp. $K^{IS}_{ij}$ or $K^{IO}_{ij}$) is non-zero, agent~$j$ should send the corresponding component of $x_j(k)$ (resp. $d_j(k)$) to agent~$i$.\
    %Notice that $E^C$, $S^x$ and $S^u$ characterize the existing communication and sensing capacities, respectively. Define matrix $\mathcal{L}\in\{0,1\}^{(p+l)\times (n+p)}$ such that $\mathcal{L}_{\ell\ell'}=0$ if there exist both a communication link and a sensor at position $(\ell,\ell')$ in the original system and $\mathcal{L}_{\ell\ell'}=1$ otherwise. If $\mathcal{L}_{\ell\ell'}=0$, adding perturbation at position $(\ell,\ell')$ generates zero additional cost. If $\mathcal{L}_{\ell\ell'}=1$, adding perturbation at position $(\ell,\ell')$ generates one additional cost (unit cost). 
    Minimizing such cost can be encoded into maximizing the sparsity of the perturbation matrix $K$ and equivalently minimizing the $\ell_0$ norm of $K$, i.e., $c(K)=\|K\|_{0}$ in problem \eqref{OP agent-wise}. This is referred to as the $\ell_{0}$ minimization and denoted by $\mathbb{P}_{0}$.

Problem $\mathbb{P}_2$: utility--privacy tradeoff. The goal of the data requester is to collect the true output. In this paper setup, the true output $y(k)$ is the linear combination of $x(k)$ and $u(k)$ weighted by the original output matrices $(\bar G,\bar H)$, i.e., $\bar Gx(k) + \bar Hu(k)$. The difference between the released output $y(k)$ of \eqref{e12} and the true output $\bar Gx(k) + \bar Hu(k)$ is data disutility. 
%So data disutility is determined by the perturbation added into the output equation \eqref{e9}. 
Notice that the perturbation added into the state equation \eqref{state eq} does not change the linear combination and thus does not affect data disutility. Instead, these perturbations can protect the privacy of target entries. We rewrite \eqref{e12} as $y(k) = \bar Gx(k)+\bar Hu(k) + [\bar H,\Pi]K[x(k)^T,u(k)^T]^T$ 
%\begin{align*}
%\left[ {\begin{array}{*{20}{c}}
%{x(k + 1)}\\
%{y(k)}
%\end{array}} \right] = \left[ {\begin{array}{*{20}{c}}
%A&B\\
%G&H
%\end{array}} \right]\left[ {\begin{array}{*{20}{c}}
%{x(k)}\\
%{y(k)}
%\end{array}} \right] + \underbrace{{\left[ {\begin{array}{*{20}{c}}
%B&\textbf{0}_{n\times l}\\
%H&\Pi
%\end{array}} \right]K\left[ {\begin{array}{*{20}{c}}
%{x(k)}\\
%{u(k)}
%\end{array}} \right]}}_\text{disturbance}.
%\end{align*}
and define data disutility as $\left\|[\bar H,\Pi]K[x(k)^T,u(k)^T]^T\right\|_2$. Notice that $\left\|[\bar H,\Pi]K[x(k)^T,u(k)^T]^T\right\|_2\leq\left\|[\bar H,\Pi]K\right\|_2$ $\left\|[x(k)^T,u(k)^T]^T\right\|_2$, and $x(k)$ and $u(k)$ are not decision variables. Hence, we turn to minimize $\left\|[\bar H,\Pi]K\right\|_2$, i.e., $c(K)=\left\|[\bar H,\Pi]K\right\|_2$ in problem \eqref{OP agent-wise}. This is referred to as the $\ell_2$ minimization and denoted by $\mathbb{P}_2$.

We assume that the optimal perturbation matrix $K$, i.e., the solution of problem \eqref{OP agent-wise}, is known to the data requester. This is another piece of auxiliary information available to the data requester.

\begin{remark}
\label{remark tradeoff}
Similar to differential privacy, our technique also adopts perturbations for privacy preservation. However, by \eqref{e36}, it can be seen that the perturbations are added in a closed-loop fashion and diminishing as the system is stabilized. Since we formulate problem \eqref{OP agent-wise} such that the perturbed system remains controllable, one can design a feedback controller by the perturbed system matrices $(\bar A+\bar BK_{SS},\bar B(I_p+K_{SI}))$ to achieve perfect stability where the perturbations vanish at the equilibrium.

Data privacy has a fundamental utility-privacy tradeoff: disclosing fully accurate information maximizes data utility but minimizes data privacy, while disclosing random noises achieves the opposite \citep{TL-NL:2009}. Our optimization formulation \eqref{OP agent-wise} utilizes control theory to characterize the tradeoff. This allows us to take into account dynamic system utilities, e.g., controllability, which have not been addressed in the literature.\oprocend
\end{remark}

\subsection{Relaxation of problem \eqref{OP agent-wise}\label{rank relaxation section}}

The first constraint of \eqref{OP agent-wise} has a clear privacy interpretation, but is not analytically tractable. In this subsection, we identify a relation between the privacy constraint and the rank deficiency of a matrix pencil, which allows us to relax the privacy constraint by a rank constraint.
%The first constraint of problem \eqref{OP agent-wise} is not in a form of any familiar problem and in general it could be difficult to directly work on this constraint. In this subsection, we relax this constraint by a rank constraint so that we can use the tools in rank-constraint optimization and matrix rank theory to systematically study the relaxed problem. To bring up the relaxation, we first introduce a sufficient condition for Definition \ref{def: agent-wise privacy}.

%\subsubsection{A sufficient condition for Definition \ref{def: agent-wise privacy}}

Given a linear system $(A,B,G,H)$, for any $z\in\mathbb{C}$, define matrix pencil\\
$D_{A,B,G,H}(z)=\left[\begin{array}{cc}
                               zI_{n}-A & -B \\
                               G & H \\
                             \end{array}
                           \right]$.\\
%                            and is used to characterize the notion of system strong observability. The linear system \eqref{state eq} and \eqref{e9} or equivalently $\left[
%                                        \begin{array}{cc}
%                                          A & B \\
%                                          G & H \\
%                                        \end{array}
%                                      \right]$ is said to be \emph{strongly observable} if, for any time instant $k$, the initial state $x(0)$ and the input sequence up to time $k-1$, i.e., $\{u(0),\cdots,u(k-1)\}$,
%can be uniquely determined from the measured output sequence $\{y(0),\cdots,y(k)\}$. By \cite{SY-MZ-EF:2016}, the linear system $\left[\begin{array}{cc}
%A & B \\
%G & H \\
%\end{array}
%\right]$ is strongly observable if and only if ${\rm rank}(D(z))=n+p$ for any $z \in \mathbb{C}$.
For any $v\in\mathbb{R}^{n+p}$, we write $v=[v_1^T,v_2^T]^T$ with $v_1\in\mathbb{R}^n$ and $v_2\in\mathbb{R}^p$. Let $v_1^t$ (resp. $v_2^t$) be the sub-vector of $v_1$ (resp. $v_2$) corresponding to $x^t(0)$ (resp. $u^t$), i.e., if the $\ell$-th entry of $x(0)$ (resp. $u$) is an entry of $x^t(0)$ (resp. $u^t$), then the $\ell$-th entry of $v_1$ (resp. $v_2$) is an entry of $v_1^t$ (resp. $v_2^t$). The dimensions of $v_1^t$ and $v_2^t$ are then $d_x^t$ and $d_u^t$, respectively. Denote by $v_{1\ell}^t$ (resp. $v_{2\ell}^t$) the $\ell$-th entry of $v_1^t$ (resp. $v_2^t$).

The following lemma provides a sufficient condition for privacy protection. Its proof leverages properties of the matrix pencil defined above, and closely follows and extends the rank-based characterizations of strong observability \citep{MLJH:1983,WK:1995}.

\begin{lemma}
\label{corollary single data item}
Given a linear system $(A,B,G,H)$, the privacy of $x^t(0)$ and $u^t$ is protected if there exists a pair of $z\in\mathbb{C}$ and $v\in\mathbb{R}^{n+p}\backslash\{\textbf{0}_{n+p}\}$ satisfying $D_{A,B,G,H}(z)v=\textbf{0}_{n+q}$ such that the following two conditions are satisfied simultaneously:\\
(1) if $d_x^t\neq0$, then $v_{1\ell}^t\neq0$ for all $\ell\in\{1,\cdots,d_x^t\}$;\\
(2) if $d_u^t\neq0$, then $z\neq 0$ and $v_{2\ell}^t\neq0$ for all $\ell\in\{1,\cdots,d_u^t\}$.\oprocend
%
%A sufficient condition that the privacy of the $\ell$-th entry of the input $u(k)$ at any time instant $k\in\mathbb{N}$ is protected is that there exists a pair of $z\in\mathbb{C}$ and $v\in\mathbb{R}^{n+p}\backslash\{\textbf{0}_{n+q}\}$ satisfying $\left[
%\begin{array}{cc}
%zI_{n}-A & -B \\
%G & H \\
%\end{array}
%\right]v=\textbf{0}_{n+q}$ such that $z\neq0$ and the $\ell$-th entry of $v_2$ is non-zero.\oprocend
\end{lemma}

\textbf{Proof:} Given that $D_{A,B,G,H}(z)v=\textbf{0}_{n+q}$, we have $Av_1+Bv_2=zv_1$ and $Gv_1+Hv_2=0_q$. Fix any $\kappa\in\mathbb{N}$ and any feasible output sequence $y_{[0,\kappa]}$. Denote by $x(0)'$ and $u_{[0,\kappa]}{'}$ an arbitrary set of initial states and input sequence that satisfy $y_{[0,\kappa]}$, i.e., $x(k+1)'=Ax(k)'+Bu(k)'$ and $y(k)=Gx(k)'+Hu(k)'$ for any $k\in\{0,\cdots,\kappa\}$. We then have $\{x^t(0)',u^{t}_{[0,\kappa]}{'}\}\in\Delta_{A,B,G,H}(y_{[0,\kappa]})$. Denote $x(0)''=x(0)'+mv_1$ and $u(k)''=u(k)'+mz^kv_2$ for each $k\in\{0,\cdots,\kappa\}$, where $m$ is an arbitrary scalar. We next show by mathematical induction that, with the initial state $x(0)''$ and input sequence $u_{[0,\kappa]}{''}$, $x(k)''=x(k)'+mz^kv_1$ for any $k\in\{0,\cdots,\kappa\}$. For $k=0$, we have $x(0)''=x(0)'+mv_1=x(0)'+mz^0v_1$. For $k=1$, we have\\
$x(1)''=Ax(0)''+Bu(0)''\\
    =A(x(0)'+mv_1)+B( u(0)'+mv_2)\\
    =Ax(0)'+Bu(0)'+m(Av_1+Bv_2)=x(1)'+mzv_1$.
    Assume that $x(k)''=x(k)'+mz^kv_1$. Then, we have\\
    $x(k+1)''=Ax(k)''+Bu(k)''\\
    =A( x(k)'+mz^kv_1)+B(u(k)'+mz^kv_2)\\
    =Ax(k)'+Bu(k)'+mz^k(Av_1+Bv_2)\\
    =x(k+1)'+mz^{k+1}v_1$.\\
    We then have $x(k)''=x(k)'+mz^kv_1$ for any $k\in\{0,\cdots,\kappa\}$. Hence, for any $k\in\{0,\cdots,\kappa\}$, we have\\
    $Gx(k)''+Hu(k)''\\
    =G(x(k)'+mz^kv_1)+H(u(k)'+mz^kv_2)\\
    =Gx(k)'+H u(k)'+mz^k(Gv_1+Hv_2)\\
    =Gx(k)'+Hu(k)'=y(k)$.\\
    This implies $\{x^t(0)'',u^t_{[0,\kappa]}{''}\}\in\Delta_{A,B,G,H}(y_{[0,\kappa]})$. Note\\
    $\min_{\ell\in\{1,\cdots,d_x^t\}} |x_{\ell}^t(0)'-x_\ell^t(0)''|=\min_{\ell\in\{1,\cdots,d_x^t\}}|mv_{1\ell}^t|,\\
    \min_{\scriptstyle \ell\in\{1,\cdots,d_u^t\} \hfill\atop \scriptstyle k\in\{0,\cdots,\kappa\} \hfill} |u_\ell^t(k)'-u_\ell^t(k)''|= \min_{\scriptstyle \ell\in\{1,\cdots,d_u^t\} \hfill\atop \scriptstyle k\in\{0,\cdots,\kappa\} \hfill}|mz^kv_{2\ell}^t|$.\\
    If $v_{1\ell}^t\neq0$ $\forall\ell\in\{1,\cdots,d_x^t\}$ and $v_{2\ell}^t\neq0$ $\forall\ell\in\{1,\cdots,d_u^t\}$ and $z\neq0$, then\\
    ${\rm Diam}_{A,B,G,H}(y_{[0,\kappa]})\geq\\
    \sup\limits_{m\in\mathbb{R}}\min\{\min\limits_{\ell\in\{1,\cdots,d_x^t\}}|mv_{1\ell}^t|,\min\limits_{\scriptstyle \ell\in\{1,\cdots,d_u^t\} \hfill\atop \scriptstyle k\in\{0,\cdots,\kappa\} \hfill}|mz^kv_{2\ell}^t|\}=\infty$.
The above analysis holds for any $\kappa\in\mathbb{N}$ and any feasible $y_{[0,\kappa]}$. By Definition \ref{def: agent-wise privacy}, $x^t(0)$ and $u^t$ are protected.\oprocend

%\subsubsection{Relaxation by rank deficiency\label{rank deficiency subsubsection}}

Lemma \ref{corollary single data item} requires that the matrix pencil $D_{A,B,G,H}(z)$ does not have full column rank. Intuitively, one can protect more entries of $x(0)$ and $u$ by reducing the rank of $D_{A,B,G,H}(z)$. This is verified by the following lemma.

\begin{lemma}
\label{rank corollary}
Given $(A,B,G,H)$, if there exists $z\neq 0$ such that $D_{A,B,G,H}(z)$ has column rank $r$, then at least $n+p-r$ entries of $x(0)$ and $u$ can be protected.\oprocend
\end{lemma}

\textbf{Proof:} If $D_{A,B,G,H}(z)$ has column rank $r$, then the null space of $D_{A,B,G,H}(z)$ has rank $n+p-r$. This implies that $D_{A,B,G,H}(z)$ must have a null vector $v$ with at least $n+p-r$ non-zero entries. By Lemma \ref{corollary single data item}, at least $n+p-r$ entries of $x(0)$ and $u$ are protected.
\oprocend

For convenience of notation, in the rest of the paper, let $\bar D(z)=D_{\bar A,\bar B,\bar G,\bar H}(z)$ and $\hat D(z)=D_{\hat A,\hat B,\hat G,\hat H}(z)$, i.e., $\bar D(z)$ and $\hat D(z)$ are the matrix pencils of the original system $(\bar A,\bar B,\bar G,\bar H)$ and the perturbed system $(\hat A,\hat B,\hat G,\hat H)$, respectively. It can be checked that $\hat D(z)=\bar D(z)+FK$ with $F = \left[\begin{array}{cc}
                               -\bar B & \textbf{0}_{n\times l} \\
                               \bar H & \Pi \\
                             \end{array}
                           \right]$. Lemma \ref{rank corollary} states that one can protect more entries of $x(0)$ and $u$ by reducing the rank of $\hat D(z)$. With more entries of $x(0)$ and $u$ being protected, in general, it is more likely that more entries of $x^t(0)$ and $u^t$ can be protected. By this observation, we relax problem \eqref{OP agent-wise} as follows:
\begin{align}
\label{relaxed OP agent-wise}
&\min\nolimits_{K\in\mathbb{R}^{(p+l)\times(n+p)},z\in\mathbb{C}}c(K)\nonumber\\
&{\rm s.t.}\;{\rm rank}(\bar D(z)+FK)<\rho,\nonumber\\
&\quad\;\; {\rm det}({\mathcal{C}}(K_{SS},K_{SI})({\mathcal{C}}(K_{SS},K_{SI}))^T) >0
%&\quad\;\; \left[\begin{array}{cc}
%                               A+B K_{SS} & B(I_{p}+K_{SI})
%                             \end{array}
%                           \right]\;\;{\rm is\;\; controllable}
\end{align}
where $\rho\in[1,\min\{n+p,n+q\}]$ is a constant integer. In the remaining, we use $\tilde{\mathbb{P}}_{0}$ (resp. $\tilde{\mathbb{P}}_2$) to denote problem \eqref{relaxed OP agent-wise} with $c(K)=\|K\|_{0}$ (resp. $c(K)=\left\|[\bar H,\Pi]K\right\|_2$).

Given any integer $\rho$ between $[1,\min\{n+p,n+q\}]$, by Lemma \ref{rank corollary}, the optimal solution of problem \eqref{relaxed OP agent-wise} can guarantee that at least $n+p-\rho+1$ entries of $x(0)$ and $u$ can be protected in the perturbed system. However, Lemma \ref{rank corollary} does not indicate which entries of $x(0)$ and $u$ can be protected. We will provide a scheme to perform the verification in the next paragraph. If some entries of $x^t(0)$ and $u^t$ are not protected, we decrease the value of $\rho$ and re-solve problem \eqref{relaxed OP agent-wise}. Our objective is to protect all the entries of $x^t(0)$ and $u^t$ with the largest possible $\rho$ (so that with the smallest possible perturbation).

We next illustrate a mechanism for checking which entries can be protected in the perturbed system after solving problem \eqref{relaxed OP agent-wise} under a given $\rho$. Given the system matrices $(\bar A,\bar B,\bar G,\bar H)$ of the original system and an optimal solution $(K,z)$ of problem \eqref{relaxed OP agent-wise},
%one can use the algorithms in \cite{AEN-PVD:1982,PM-PVD-AV:1994}~\footnote{The algorithm of \cite{AEN-PVD:1982,PM-PVD-AV:1994} for finding the invariant zeros of a linear system has been implemented in Matlab by the command $tzero$ \cite{tzero}.} to obtain the invariant zeros $z$'s of the perturbed system, i.e., those $z$'s such that ${\rm rank}(D(z)+FK)<n+p$. For each invariant zero $z$,
one can derive the null space of $\bar D(z)+FK$ and then make use of Lemma \ref{corollary single data item} to check which entries can be protected in the perturbed system. In particular, if the null space of $\bar D(z)+FK$ admits a null vector $v=[v_1^T,v_2^T]^T$ such that $v_{1\ell}^t\neq0$, then the $\ell$-th entry of $x^t(0)$ is protected. If $z\neq0$ and the null space of $\bar D(z)+FK$ admits a null vector $v$ such that $v_{2\ell}^t\neq0$, then the $\ell$-th entry of $u^t$ is protected.
%To check whether the $\ell$-th entry of $x(0)$ is protected, one can check whether the null space of $D(z)+FK$ admits a null vector $v=[v_1^T,v_2^T]^T$ such that the $\ell$-th entry of $v_1$ is non-zero. To check whether the $\ell$-th entry of $u(k)$ is protected, one can check whether there exists a non-zero invariant zero $z$ such that the corresponding null space of $D(z)+FK$ admits a null vector $v$ such that the $\ell$-th entry of $v_2$ is non-zero.

The remaining issue is how to solve problem \eqref{relaxed OP agent-wise} under a given $\rho$. The next theorem shows the non-convexity of the constraint set of problem \eqref{relaxed OP agent-wise}, indicating that the problem could be hard to solve and needs to be further relaxed. In the next section, we furhter prove a NP-hardness result for problem $\tilde{\mathbb{P}}_{0}$.

\begin{theorem}
\label{theorem: nonconvex constraint}
%Assume $\mathbb{K}=\mathbb{R}^{(p+l)\times(n+p)}$. 
The constraint set of problem \eqref{relaxed OP agent-wise} is non-convex.\hfill\oprocend
\end{theorem}

\textbf{Proof:} Denote the constraint set of \eqref{relaxed OP agent-wise} by $\Upsilon$.
%For the controllability constraint, notice that $(A+BK_{SS},B(I_p+K_{SI}))$ is controllable if and only if the controllability matrix ${\mathcal{C}}(K_{SS},K_{SI})$ has full row rank.
Let $(z',K')\in \Upsilon$ be any feasible quadruple. The feasibility implies that ${\mathcal{C}}(K_{SS}',K_{SI}')$ has full row rank and ${\rm rank}(\bar D(z')+FK')<\rho$.
%${\rm rank}(\left[\begin{array}{cc}
%                               \bar zI_{n}-(A+B \bar K_{SS}) & -B (I_{p}+\bar K_{SI}) \\
%                               \Pi(G'+H' \bar K_{SS} + \bar K_{OS}) & \Pi(H'+H' \bar K_{SI} + \bar K_{OI}) \\
%                             \end{array}
%                           \right])<\rho$.
Consider $(K_{SS}'',K_{SI}'',K_{OS}'',K_{OI}'')$ $=(K_{SS}',-K_{SI}'-2I_{p},K_{OS}',-K_{OI}')$. Then\\
${\mathcal{C}}(K_{SS}'', K_{SI}'')={\mathcal{C}}(K_{SS}',-K_{SI}'-2I_{p}) \\
= [\bar B(I_{p}-K_{SI}'-2I_{p}),(\bar A+\bar BK_{SS}')\bar B(I_{p}-K_{SI}'-2I_{p}),\\
\cdots,(\bar A+\bar B K_{SS}')^{n-1}\bar B(I_{p}-K_{SI}'-2I_{p})]\\
=-{\mathcal{C}}(K_{SS}',K_{SI}')$\\
which has full row rank. Thus, $(K_{SS}'',$ $K_{SI}'')$ satisfies the controllability constraint. Take $z''=z'$. We have\\
$\bar D(z'')+FK''\\
=\left[\begin{array}{cc}
z''I_{n}-(\bar A+\bar B K_{SS}'') & -\bar B (I_{p}+K_{SI}'') \\
\bar G+\bar H K_{SS}'' + \Pi K_{OS}'' & \bar H+\bar H K_{SI}'' + \Pi K_{OI}'' \\
\end{array}
\right]\\
=\left[\begin{array}{cc}
z'I_{n}-(\bar A+\bar B K_{SS}') & \bar B (I_{p}+K_{SI}') \\
\bar G+\bar H K_{SS}' + \Pi K_{OS}' & -(\bar H+\bar H K_{SI}' + \Pi K_{OI}') \\
\end{array}
\right]$\\
which implies ${\rm rank}(\bar D(z'')+FK'')={\rm rank}(\bar D(z')+FK')<\rho$.
%clearly has the same column rank as
%\begin{align*}
%\left[\begin{array}{cc}
%\bar zI_{n}-(A+B \bar K_{SS}) & -B (I_{p}+\bar K_{SI}) \\
%\Pi(G'+H' \bar K_{SS} + \bar K_{OS}) & \Pi(H'+H' \bar K_{SI} + \bar K_{OI}) \\
%\end{array}
%\right]
%\end{align*}
%and does not have full column rank. Thus, the quadruple $(\hat K_{SS},\hat K_{SI},\hat K_{OS},\hat K_{OI})$ satisfies the non-strong-observability constraint.
Thus, $(z'',K'')\in\Upsilon$.

Now consider another quadruple $(K_{SS}^\circ,K_{SI}^\circ,K_{OS}^\circ,K_{OI}^\circ)$ $=(K_{SS}',-I_{p},K_{OS}',\textbf{0}_{l\times p})$. Notice that $(K_{SS}^\circ,K_{SI}^\circ,K_{OS}^\circ,$ $K_{OI}^\circ)=(K_{SS}',-I_{p},K_{OS}',\textbf{0}_{l\times p})=\frac{1}{2}(K_{SS}',K_{SI}',K_{OS}',$ $K_{OI}')+\frac{1}{2}(K_{SS}'',K_{SI}'',K_{OS}'',K_{OI}'')$, which implies that $(K_{SS}^\circ,K_{SI}^\circ,K_{OS}^\circ,K_{OI}^\circ)$ is a convex combination of $(K_{SS}',K_{SI}',K_{OS}',K_{OI}')$ and $(K_{SS}'',K_{SI}'',$ $K_{OS}'',K_{OI}'')$. Let $z^\circ=z'$. Then $(z^\circ,K^\circ)$ is a convex combination of $(z',K')$ and $( z'',K'')$. We have\\
${\mathcal{C}}(K_{SS}^\circ,K_{SI}^\circ)={\mathcal{C}}(K_{SS}',-I_{p}) \\
= [\bar B(I_{p}-I_{p}),(\bar A+\bar B K_{SS}')\bar B(I_{p}-I_{p}),\cdots,\\
\quad\;(\bar A+\bar BK_{SS}')^{n-1}\bar B(I_{p}-I_{p})]=\textbf{0}_{n\times np}$,\\
which does not have full row rank. Thus, $(K_{SS}^\circ,K_{SI}^\circ)$ does not satisfy the controllability constraint and hence $(z^\circ,K^\circ)\notin\Upsilon$. This implies that $\Upsilon$ is non-convex.\hfill\oprocend

\section{Problem $\tilde{\mathbb{P}}_{0}$\label{P0 section}}

In this section, we first prove that a relaxation of problem $\tilde{\mathbb{P}}_{0}$ is NP-hard, which indicates that problem $\tilde{\mathbb{P}}_{0}$ itself might also be NP-hard. We then provide a convex approximation for problem $\tilde{\mathbb{P}}_{0}$. Specifically, in problem \eqref{relaxed OP agent-wise}, the $\ell_0$ norm in the objective function is relaxed by the $\ell_1$ norm heuristic, the rank constraint is relaxed by the nuclear norm heuristic, and the controllability constraint is approximated by a symmetric positive semidefinite condition.

\subsection{Computational intractability}

To obtain a rigorous NP-hardness result, we consider the following problem derived by fixing some $z$ and $[K_{SS},K_{SI}]=\textbf{0}_{p\times(n+p)}$, and dropping the controllability constraint of problem $\tilde{\mathbb{P}}_{0}$:
%can be equivalently written in the following bilevel optimization structure:
%\begin{align}
%\label{relaxed OP bilevel}
%&\min\nolimits_{z\in\mathbb{C}}\tilde{\mathbb{P}}_{0,\mathcal{L}}(z)\nonumber\\
%&{\rm s.t.}\;\tilde{\mathbb{P}}_{0,\mathcal{L}}(z)=\mathop {\min }\limits_{\begin{array}{*{20}c}
%   \substack{\{K\in\mathbb{K}:{\rm rank}(\bar D(z)+FK)<\rho,\\{\rm det}({\mathcal{C}}(K_{SS},K_{SI})({\mathcal{C}}(K_{SS},K_{SI}))^T) >0\}}\\
% \end{array} }\|K\|_{0,\mathcal{L}}.
%%&\quad\;\; \left[\begin{array}{cc}
%%                               A+B K_{SS} & B(I_{p}+K_{SI})
%%                             \end{array}
%%                           \right]\;\;{\rm is\;\; controllable}
%\end{align}
%Consider the following problem derived by dropping the controllability constraint of problem $\tilde{\mathbb{P}}_{0,\mathcal{L}}(z)$:
\begin{align}
\label{special lower level OP}
&\min\nolimits_{K_O\in\mathbb{R}^{l\times(n+p)}}\|K_O\|_{0}\nonumber\\
&{\rm s.t.}\;\;{\rm rank}(\bar D(z)+F\left[ {\begin{array}{*{20}{c}}
{{\textbf{0}_{p \times (n + p)}}}\\
{{K_O}}
\end{array}} \right])<\rho.
\end{align}
where $K_O=[K_{OS},K_{OI}]$. Denote problem \eqref{special lower level OP} by $\hat{\mathbb{P}}_{0}(z)$. By fixing $z$ and $[K_{SS},K_{SI}]$, the dimension of the decision variables is reduced. By the proof of Theorem \ref{theorem: nonconvex constraint}, the controllability constraint of problem $\tilde{\mathbb{P}}_{0}$ is non-convex. Hence, intuitively, problem $\tilde{\mathbb{P}}_{0}$ might be harder to solve than problem $\hat{\mathbb{P}}_{0}(z)$. We next show that problem $\hat{\mathbb{P}}_{0}(z)$ is NP-hard due to the non-convexity of its objective function. This provides an implication that problem $\tilde{\mathbb{P}}_{0}$ might also be NP-hard. We leave the proof of NP-hardness of problem $\tilde{\mathbb{P}}_{0}$ to our future works.
%Moreover, in general, a bilevel optimization problem could be harder to solve than its lower level problem. For example, bilevel linear programming is NP-hard despite the lower level linear programming can be solved efficiently \cite{RGJ:1985}. Hence, the NP-hardness of problem $\hat{\mathbb{P}}_{0,\mathcal{L}}(z)$ provides an implication that problem $\tilde{\mathbb{P}}_{0,\mathcal{L}}$ could be NP-hard.

%We next show that the lower level problem $\tilde{\mathbb{P}}_{0,\mathcal{L}}(z)$ is NP-hard due to the non-convexity of its objective function.
It is well-known that $\ell_0$ norm is non-convex and $\ell_0$ norm optimization problems are hard to solve in general. However, there has been a limited number of $\ell_0$ norm optimization problems which have been rigorously proven to be NP-hard. The following theorem establishes the NP-hardness of problem $\hat{\mathbb{P}}_{0}(z)$ by showing that it is as hard as finding a sparsest null vector of a matrix with more columns than rows, which has been proven to be NP-hard \citep{TC-AP:1986}.

\begin{theorem}
Problem $\hat{\mathbb{P}}_{0}(z)$ is NP-hard.\hfill\oprocend
\label{theorem: NP hard}
\end{theorem}
%
%\begin{lemma}
%Problem $\mathbb{P}_0$ is NP-hard.
%\end{lemma}
%
%\begin{assumption}
%$[A\;\;G]$ is controllable. \hfill\oprocend
%\label{asm: controllable}
%\end{assumption}

\textbf{Proof:} We prove the NP-hardness of problem $\hat{\mathbb{P}}_{0}(z)$ by following the standard procedure of proving NP-hardness \citep{JVL:1990}:\\
\emph{Step 1}. Reduce any instance of a known NP-hard problem to an instance of problem $\hat{\mathbb{P}}_{0}(z)$ in polynomial time;\\
\emph{Step 2}. Show that a solution of the instance of the known NP-hard problem can be constructed from a solution of the instance of problem $\hat{\mathbb{P}}_{0}(z)$ in polynomial time.

%We prove the NP-hardness of problem $\hat{\mathbb{P}}_{0,\mathcal{L}}(z)$ by reducing in polynomial time a known NP-hard problem to it.
The known NP-hard problem we use is the following null vector problem (NVP) \citep{TC-AP:1986}:

\textbf{Null vector problem:} Given a matrix $M\in\mathbb{R}^{r\times c}$ with $r<c$, find a sparest null vector of $M$, i.e., find an optimal solution $v^*$ to the following optimization problem
\begin{flalign}
&\min\nolimits_{v\in\mathbb{R}^c\backslash\{\textbf{0}_c\}}\|v\|_0\quad{\rm s.t.}\;\; Mv=\textbf{0}_{r}.
\label{e3_3}
\end{flalign}
%$v^*={\rm argmin}_{v\in\{\mathbb{R}^c\backslash\textbf{0}_c:Mv=\textbf{0}_r\}}\|v\|_0$.
%By the abstract of \cite{TC-AP:1986}, the NVP is NP-hard.
%We next show that the NVP can be reduced in polynomial time to problem $\hat{\mathbb{P}}_{0,\mathcal{L}}(z)$.
\emph{Step 1}. Consider a matrix $M\in\mathbb{R}^{r\times c}$ with $r<c$. Let $M=[M_1,M_2]$ with $M_1\in\mathbb{R}^{r\times r}$ and $M_2\in\mathbb{R}^{r\times (c-r)}$. Given $M$, we construct an instance of $\hat{\mathbb{P}}_{0}(z)$ as follows: let $n=r$, $p=c-r$, $l=q=c$, $z$ be any fixed complex number, $\bar A=zI_n-M_1$, $\bar B=-M_2$, $[\bar G',\bar H']=\Pi=I_{n+p}$, 
%$\mathbb{K}$ be such that $K_{SS}=\textbf{0}_{p\times n}$, $K_{SI}=\textbf{0}_{p\times p}$, $K_{OS}\in\mathbb{R}^{l\times n}$ and $K_{OI}\in\mathbb{R}^{l\times p}$,
and $\rho=n+p$. 
%and $\mathcal{L}_{\ell\ell'}=1$ for all $(\ell,\ell')$'s corresponding to $K_{OS}$ and $K_{OI}$.
%The restriction of $K_{SS}=\textbf{0}_{p\times n}$, $K_{SI}=\textbf{0}_{p\times p}$, $K_{OS}\in\mathbb{R}^{l\times n}$ and $K_{OI}\in\mathbb{R}^{l\times p}$ can be encoded in the structure of $\mathbb{K}$.
%We next formulate problem $\hat{\mathbb{P}}_{0,\mathcal{L}}(z)$ with the above parameters.
%Since $K_{SS}$ and $K_{SI}$ are zero matrices and $(\bar A,\bar B)$ is controllable, the controllability constraint of $\tilde{\mathbb{P}}_{0,\mathcal{L}}(z)$ is satisfied and can be removed.
%Since $z=0$ and $\Pi=I_{n+p}$,
%and $\Pi$ is non-singular, we have ${\rm rank}(\left[\begin{array}{cc}
%                               zI_{n}-\bar A & -\bar B \\
%                               \bar G & \bar H \\
%                             \end{array}
%                           \right])=n+p$ for any $z\in\mathbb{C}$. Since $\Pi$ is non-singular, for any $z\in\mathbb{C}$,
%we have $\left[\begin{array}{cc}
%                               zI_{n}-\bar A & -\bar B \\
%                               \Pi(\bar G'+K_{OS}) & \Pi(\bar H'+K_{OI}) \\
%                             \end{array}
%                           \right]=\left[\begin{array}{cc}
%                               -\bar A & -\bar B \\
%                               \bar G'+K_{OS} & \bar H'+K_{OI} \\
%                             \end{array}
%                           \right]$.
%Let $K_O=[K_{OS},K_{OI}]$. 
With the above defined parameters, the instance of problem $\hat{\mathbb{P}}_{0}(z)$ can be written as:
\begin{align}
&\!\!\!\!\!\min\limits_{K_O\in\mathbb{R}^{c\times c}}\|K_O\|_0\;\,{\rm s.t.}\, {\rm rank}([M^T,(I_c+K_O)^T]^T)<c.%\nnum\\
%&\quad\quad \left[\begin{array}{cc}
%                               A & G
%                             \end{array}
%                           \right] + \left[\begin{array}{cc}
%                               G K_{SS} & G K_{SI}
%                             \end{array}
%                           \right]\nnum\\
%&\quad\quad{\rm is\;\;controllable}
\label{e3_2}
\end{align}
%Consider another related problem cast as follows:
%\begin{flalign}
%&\min\limits_{v\in\mathbb{R}^c\backslash\{\textbf{0}_c\}}\|v\|_0\quad{\rm s.t.}\;\; Mv=\textbf{0}_{r}.
%\label{e3_3}
%\end{flalign}
It is clear that the construction of problem \eqref{e3_2} can be done in polynomial time.

\emph{Step 2}. Let $K_O^*$ be an optimal solution of problem \eqref{e3_2} and $v^*$ be an optimal solution of problem \eqref{e3_3}.
%The following claim determines a relation between the two solutions.
%We first prove the following claim.

\textbf{Claim I:} $\|K_O^*\|_0=\|v^*\|_0$.

\textbf{Proof of Claim I:} Let $K_O$ be any feasible solution of problem \eqref{e3_2}. Notice that ${\rm rank}([M^T,(I_c+K_O)^T]^T)<c$ if and only if there exists $v\in\mathbb{R}^c\backslash\{\textbf{0}_c\}$ such that $Mv=\textbf{0}_r$ and $K_Ov=-v$. For any such vector $v$, for any $i\in\{1,\cdots,c\}$, if $v_i\neq0$, then, to satisfy $K_Ov=-v$, the entries of the $i$-th row of $K_O$ cannot be all zero. This implies $\|K_O\|_0\geq\|v\|_0$. In particular, since $v$ is a null vector of $M$, we have $\|K_O\|_0\geq\|v^*\|_0$. Since this is true for any feasible $K_O$, we have $\|K_O^*\|_0\geq\|v^*\|_0$. Next, by the following procedure, we construct a matrix $\tilde K_O$ that is feasible to problem \eqref{e3_2} and $\|\tilde K_O\|_0=\|v^*\|_0$.

\emph{Procedure I}: For each $i\in\{1,\cdots,c\}$, if $v_i^*=0$, then $(\tilde K_O)_{ij}=0$ for all $j\in\{1,\cdots,c\}$; if $v_i^*\neq0$, then $(\tilde K_O)_{ii}=-1$ and $(\tilde K_O)_{ij}=0$ for all $j\in\{1,\cdots,c\}\backslash\{i\}$.

By Procedure I, it is easy to derive $\tilde K_Ov^*=-v^*$ and $\|\tilde K_O\|_0=\|v^*\|_0$. Since $v^*$ is a null vector of $M$, we have $Mv^*=\textbf{0}_r$. By optimality, we have $\|K_O^*\|_0\leq\|\tilde K_O\|_0=\|v^*\|_0$. Hence, together with the above result $\|K_O^*\|_0\geq\|v^*\|_0$, we have $\|K_O^*\|_0=\|v^*\|_0$.\hfill\oprocend

%\textbf{Claim II:} An optimal solution $v^*$ of problem \eqref{e3_3} can be constructed in polynomial time from an optimal solution $K_O^*$ of problem \eqref{e3_2}.

We next complete Step 2 by showing that an optimal solution $v^*$ of problem \eqref{e3_3} can be constructed in polynomial time from an optimal solution $K_O^*$ of problem \eqref{e3_2}.
%\textbf{Proof of Claim II:}
%By Claim I, if one can obtain an optimal solution $K_O^*$ of problem \eqref{e3_2} in polynomial time, then it can derive the optimal value $\|v^*\|_0$ of problem \eqref{e3_3} in polynomial time by $\|v^*\|_0=\|K_O^*\|_0$.
Let $\eta=\|K_O^*\|_0$. Notice that $\eta$ is a known constant. We next show that with $\eta$, one can derive an optimal solution $v^*$ of problem \eqref{e3_3} in polynomial time.

%Given a vector $v$ with $m$ non-zero entries, if the $i$-th entry of $v$ is non-zero, then, to satisfy $K_Ov=-v$, the entries of the $i$-th row of $K_O$ cannot be all zero. Thus, to satisfy $K_Ov=-v$, $K_O$ must at least have $m$ non-zero entries. It is clear that $K_O$ constructed by Procedure I has exactly $m$ non-zero entries and satisfies $K_Ov=-v$. This implies that Procedure I generates a $K_O$ that satisfies $K_Ov=-v$ with the least non-zero entries. Thus, $K_O^*$ must correspond to a vector $v$ with the least number of non-zero entries under the constraint $Mv=\textbf{0}_{n}$. Since $v^*$ is the optimal solution of problem \eqref{e3_3}, $K_O^*$ constructed by Procedure I via $v^*$ has the least non-zero entries and is an optimal solution of problem \eqref{e3_2}. Hence, $\|K_O^*\|_0=\|v^*\|_0$.

By Claim I, we have $\|v^*\|_0=\eta$, i.e., a sparest null vector of $M$ has $\eta$ non-zero entries. To find a sparest null vector of $M$, we consider the sub-matrices composed of any collection of $\eta$ columns of $M$.
%Hence, any combination of $\eta-1$ or less columns of $M$ is linearly independent and there exists at least one combination of $\eta$ columns of $M$ such that the $\eta$ columns are linearly dependent. Hence, to find a sparest null vector of $M$, one only needs to enumerate all possible sub-matrices composed of the combinations of $\eta$ columns of $M$.
There are totally $\left( {\begin{array}{*{20}{c}}
{c}\\
\eta
\end{array}} \right)=\frac{(c-\eta+1)\cdot\cdots\cdot c}{\eta!}$
%$\left( {\begin{array}{*{20}{c}}
%{c}\\
%\eta
%\end{array}} \right)=\frac{c!}{\eta!(c-\eta)!}=\frac{(c-\eta+1)\cdot\cdots\cdot c}{\eta!}$
collections of $\eta$ columns, which is a polynomial of $c$. For each sub-matrix $M_\eta$ of $\eta$ columns, we solve $M_\eta v=\textbf{0}_r$ to obtain the general form of solution of $v$, which can be done in polynomial time by the Gaussian elimination method (page 12 of \citep{RWF:1988}). If the general form of solution only admits $\textbf{0}_\eta$, then go on with the next sub-matrix; if the the general form of solution admits a vector with at least one non-zero entry, stop. Denote the matrix at the last step by $M_\eta^*$. Pick any non-zero vector from the general form of solution to $M_\eta^* v=\textbf{0}_r$, and denote it by $\bar v^*$. Since $\bar v^*\in\mathbb{R}^\eta\backslash \{\textbf{0}_\eta\}$, we have $0<\|\bar v^*\|_0\leq\eta$. We then augment $\bar v^*$ to a vector $\hat v^*\in\mathbb{R}^c$ by filling in zeros to the positions corresponding to the columns of $M$ that are not in $M_\eta^*$.
%That is, if a column $M_i$ is in $M_\eta^*$, then $v_i^*$ is equal to the corresponding entry of $\bar v^*$; if $M_i$ is not in $M_\eta^*$, then $v_i^*=0$.
It is clear that $M\hat v^*=\textbf{0}_r$ and $\|\hat v^*\|_0=\|\bar v^*\|_0$. Hence, $\hat v^*$ is a null vector of $M$ such that $0<\|\hat v^*\|_0\leq\eta$. Since a sparsest null vector of $M$ has $\eta$ non-zero entries, it must hold $\|\hat v^*\|_0\geq\eta$. Then $\|\hat v^*\|_0=\eta$.
%In other words, all the entries of $\bar v^*$ are non-zero. Otherwise, there exists a combination of $\eta-1$ columns of $M$ such that these $\eta-1$ columns are linearly dependent, which contradicts $\|v^*\|_0=\eta$.
Hence, we have derived an optimal solution of problem \eqref{e3_3}. Since the total number of sub-matrices is a polynomial of $c$ and for each sub-matrix, it takes polynomial time to do the computation, we have constructed a solution to the NVP from a solution of problem $\hat{\mathbb{P}}_{0}(z)$ in polynomial time.\hfill\oprocend

\subsection{Convex relaxations}

The $\ell_0$ norm $\|\cdot\|_{0}$ in problem \eqref{relaxed OP agent-wise} introduces non-convexity. In compressed sensing \citep{EC-TT:05,DD:06}, it is a common practice to replace $\|\cdot\|_0$ by $\|\cdot\|_1$. It is proven that the $\ell_1$ norm heuristic returns the sparsest solution under certain conditions, e.g., restricted isometry property (RIP) \citep{EC-TT:05}. Through experiments, one can see that the $\ell_1$ norm heuristic can return sparse solutions even RIP is not valid \citep{JY-YZ:2011}.
%We stack the entries $T_{\ell \ell'}$ of matrix $T=K\circ\mathcal{L}$ into a vector denoted by ${\rm vec}(K)_{\mathcal{L}}$ and denote the dimension of ${\rm vec}(K)_{\mathcal{L}}$ by $m$.
Recall that ${\rm vec}(K)$ is the column vector consisting of the entries of $K$. By the $\ell_1$ norm relaxation, the objective function $\|K\|_{0}$ is relaxed by $\|{\rm vec}(K)\|_1$.

The constraint ${\rm rank}(\bar D(z) + FK) < \rho$ is a rank constraint. In general, rank constraint/minimization problems is hard to solve, both in theory and practice. A particularly interesting method is the nuclear norm heuristic. In particular, \citep{MF-HH-SB:2004} showed that the convex envelop of the function ${\rm rank}(M)$ on the set $\{M\in\mathbb{R}^{m\times n}|\|M\|_2\leq1\}$ is $\|M\|_*$. In addition, \citep{BR-MF-PAP:2010} showed that, under certain conditions, e.g., RIP, the relaxation via the nuclear norm heuristic can return minimum-rank solutions. By the nuclear norm relaxation, the rank constraint of problem \eqref{relaxed OP agent-wise} is relaxed by $\min_{z\in\mathbb{C},K\in\mathbb{R}^{(p+l)\times(n+p)}}\|\bar D(z)+FK\|_*$. Notice that this relaxation turns the hard rank constraint into a soft constraint in the objective function.

The determinant function in the second constraint of problem \eqref{relaxed OP agent-wise} is a polynomial of the entries of $K_{SS}$ and $K_{SI}$ and is non-convex. To relax this controllability constraint, we first introduce the following lemma.

\begin{lemma}
Assume that $(\bar A,\bar B)$ is controllable. Then $(\hat A,\hat B)$ is controllable if and only if $v\bar BK_{SI}\neq-v\bar B$ for any left eigenvector $v$ of $\bar A+\bar BK_{SS}$.\hfill\oprocend
%i.e., $v\in\mathbb{C}^{1\times n}$ and $v\neq\textbf{0}_{1\times n}$ and there exists $\lambda\in\mathbb{C}$ such that $v(\bar A+\bar BK_{SS})=\lambda (\bar A+\bar BK_{SS})$.
\label{lemma: controllale}
\end{lemma}

\textbf{Proof:} %By Theorem 8.M1 of \cite{Chen}, the perturbed system $(A+BK_{SS},B(I_p+K_{SI}))$ is controllable if and only if the system $(A,B(I_p+K_{SI}))$ is controllable.
By Theorem 6.1 of \citep{Chen}, the perturbed system $(\bar A+\bar BK_{SS},\bar B(I_p+K_{SI}))$ is controllable if and only if $[\bar A+\bar BK_{SS} - \lambda I_n ,\bar B(I_p+K_{SI})]$ has full row rank at every eigenvalue $\lambda$ of $\bar A+\bar BK_{SS}$, or, equivalently, $v[\bar A+\bar BK_{SS} - \lambda I_n ,\bar B(I_p+K_{SI})]\neq\textbf{0}_{1\times (n+p)}$ for any $v\in\mathbb{C}^{1\times n}$ and $v\neq\textbf{0}_{1\times n}$ at every eigenvalue $\lambda$ of $\bar A+\bar BK_{SS}$. Since $\bar A+\bar BK_{SS}$ is real and $\lambda$ is an eigenvalue of $\bar A+\bar BK_{SS}$, the latter condition above is then equivalent to the condition that for each eigenvalue $\lambda$ of $\bar A+\bar BK_{SS}$, $v[\bar A+\bar BK_{SS} - \lambda I_n ,\bar B(I_p+K_{SI})]\neq\textbf{0}_{1\times (n+p)}$ for any left eigenvector $v$ of $\bar A+\bar BK_{SS}$ corresponding to $\lambda$, i.e., $v(\bar A+\bar BK_{SS})=\lambda v$, for otherwise, if $v$ is not a left eigenvector of $\bar A+\bar BK_{SS}$ corresponding to $\lambda$, then $v(\bar A+\bar BK_{SS} - \lambda I_n)\neq\textbf{0}_{1\times n}$ and surely $v[\bar A+\bar BK_{SS} - \lambda I_n ,\bar B(I_p+K_{SI})]\neq\textbf{0}_{1\times (n+p)}$. For each eigenvalue $\lambda$ of $\bar A+\bar BK_{SS}$, the condition $v[\bar A+\bar BK_{SS} - \lambda I_n ,\bar B(I_p+K_{SI})]\neq\textbf{0}_{1\times (n+p)}$ for any left eigenvector $v$ of $\bar A+\bar BK_{SS}$ corresponding to $\lambda$ is equivalent to $v\bar B(I_p+K_{SI})\neq\textbf{0}_{1\times p}$, and further equivalent to $v\bar BK_{SI}\neq-v\bar B$ for any left eigenvector $v$ of $\bar A+\bar BK_{SS}$ corresponding to $\lambda$. Since this needs to hold for every eigenvalue $\lambda$ of $\bar A+\bar BK_{SS}$, we have that $(\hat A,\hat B)$ is controllable if and only if $v\bar BK_{SI}\neq-v\bar B$ for any left eigenvector $v$ of $\bar A+\bar BK_{SS}$.\hfill\oprocend

%The latter condition is possible if and only if $v\bar B\neq\textbf{0}_{1\times p}$ for any left eigenvector $v$ of $\bar A+\bar BK_{SS}$. We next show that this is true. Since $(\bar A,\bar B)$ is controllable, by Theorem 8.M1 of \cite{Chen}, $(\bar A+\bar BK_{SS},\bar B)$ is also controllable. Thus, $\left[ {\begin{array}{*{20}c}
%   {\bar A+\bar BK_{SS} - \lambda I_n } & \bar B  \\
% \end{array} } \right]$ has full row rank at every eigenvalue $\lambda$ of $\bar A+\bar BK_{SS}$, which implies that $v\bar B\neq\textbf{0}_{1\times p}$ for any left eigenvector $v$ of $\bar A+\bar BK_{SS}$.

%If $I_p+K_{SI}$ is invertible, then $w(I_p+K_{SI})\neq\textbf{0}_{1\times p}$ for any $w\in\mathbb{R}^{1\times p}\backslash\{\textbf{0}_{1\times p}\}$.
Corollary \ref{corollary controllability} states that the invertibility of $I_p+K_{SI}$ is a sufficient condition for the controllability of $(\hat A,\hat B)$.

\begin{corollary}
Assume that $(\bar A,\bar B)$ is controllable. If $I_p+K_{SI}$ is invertible, then $(\hat A,\hat B)$ is controllable.\hfill\oprocend
\label{corollary controllability}
\end{corollary}

\textbf{Proof:} Since $(\bar A,\bar B)$ is controllable, by Theorem 8.M1 of \citep{Chen}, $(\bar A+\bar BK_{SS},\bar B)$ is controllable. Thus, $[\bar A+\bar BK_{SS} - \lambda I_n,\bar B]$ has full row rank at every eigenvalue $\lambda$ of $\bar A+\bar BK_{SS}$. So $v\bar B\neq\textbf{0}_{1\times p}$ for any left eigenvector $v$ of $\bar A+\bar BK_{SS}$. Hence, if $I_p+K_{SI}$ is invertible, then $v\bar B(I_p+K_{SI})\neq\textbf{0}_{1\times p}$ for any left eigenvector $v$ of $\bar A+\bar BK_{SS}$. By Lemma \ref{lemma: controllale}, we have that $(\hat A,\hat B)$ is controllable.\hfill\oprocend

The invertibility of $I_p+K_{SI}$ is equivalent to that its determinant is non-zero. However, the determinant of $I_p+K_{SI}$ is a polynomial of the entries of $K_{SI}$ and is non-convex. We further relax the invertibility of $I_p+K_{SI}$ by the condition that $K_{SI}$ is symmetric and $I_p+K_{SI}\succ0$. The strict positive definite condition is usually difficult to deal with and may lead to infeasibility of the problem. We relax this by a semidefinite condition as $I_p+K_{SI}-\varepsilon I_p\succeq0$, where $\varepsilon>0$ is a tuning parameter. It is clear that if $(1-\varepsilon)I_p+K_{SI}\succeq0$, then $I_p+K_{SI}$ is invertible.

With the above relaxations, $\tilde{\mathbb{P}}_{0}$ is approximated by:
\begin{align}
\label{relaxation 1}
&\min\limits_{z\in \mathbb{C},K\in\mathbb{R}^{(p+l)\times(n+p)},K_{SI}\in\mathbb{S}^{p}}\|{\rm vec}(K)\|_{1}+c\|\bar D(z) + FK\|_*\nonumber\\
&{\rm s.t.}\;\;(1-\varepsilon)I_p+K_{SI}\succeq0.
\end{align}
%where $c>0$ and $\varepsilon>0$ are tuning parameters.

\begin{remark}
\label{remark on c}
In problem \eqref{relaxation 1}, $c>0$ plays the role of Lagrangian multiplier, and tunes the relative weights between $\|{\rm vec}(K)\|_{1}$ and $\|\bar D(z) + FK\|_*$.
%Intuitively, if $c$ is large, then the minimization is biased on $\|\bar D(z) + FK\|_*$ and the perturbed system tends to lose more ranks; if $c$ is small, then the minimization is biased on $\|{\rm vec}(K\circ\mathcal{L})\|_{1}$ and the perturbation matrix $K$ tends to be sparse. This tradeoff is observed in the case study in Section \ref{simulation section}. Hence, for the relaxed problem \eqref{relaxation 1}, the tuning of $\rho$ mentioned in Section \ref{rank relaxation section} becomes the tuning of $c$.
\hfill\oprocend
\end{remark}

With the linear program (LP) characterization of $\ell_1$ norm given in page 294 of \citep{Boyd.Vandenberghe:04}, $\min_{K\in\mathbb{R}^{(p+l)\times(n+p)}} \|{\rm vec}(K)\|_1$ can be cast as: $\min\nolimits_{K\in\mathbb{R}^{(p+l)\times(n+p)},t\in\mathbb{R}^m}\sum\nolimits_{\ell=1}^mt_\ell$, ${\rm s.t.} -t\leq {\rm vec}(K)\leq t$. With the semidefinite program (SDP) characterization of nuclear norm given by \citep{BR-MF-PAP:2010}, $\min_{z\in\mathbb{C},K\in\mathbb{R}^{(p+l)\times(n+p)}} c\|\bar D(z) + FK\|_*$ can be cast as:
\begin{align*}
&\min\limits_{z\in\mathbb{C},K\in\mathbb{R}^{(p+l)\times(n+p)},W_1\in\mathbb{S}^{n+l},W_2\in\mathbb{S}^{n+p}}c({\rm Tr}(W_1)+{\rm Tr}(W_2))\nonumber\\
&{\rm s.t.}\;\;  \left[\begin{array}{cc}
                     W_1 & \bar D(z) + FK \\
                     (\bar D(z) + FK)^T & W_2 \\
                   \end{array}
                 \right]
\succeq0.
%&\quad\quad W_1 = W_1^T,\quad W_2 = W_2^T.
\end{align*}
With the above LP and SDP characterizations, problem \eqref{relaxation 1} can be equivalently turned into an SDP as follows:
\begin{align}
\label{relaxation 2}
&\mathop {\min }\limits_{\begin{array}{*{20}c}
   \substack{ z\in\mathbb{C},K\in\mathbb{R}^{(p+l)\times(n+p)},t\in\mathbb{R}^m \\ K_{SI}\in\mathbb{S}^{p},W_1  \in \mathbb{S}^{n + l}, W_2  \in \mathbb{S}^{n + p} }   \\
 \end{array} }\sum\nolimits_{\ell=1}^mt_\ell+c({\rm Tr}(W_1)+{\rm Tr}(W_2))\nonumber\\
&{\rm s.t.}\quad  -t\leq {\rm vec}(K)\leq t,\quad (1-\varepsilon)I_p+K_{SI}\succeq0\nonumber\\
&\quad\quad\left[\begin{array}{cc}
                     W_1 & \bar D(z) + FK \\
                     (\bar D(z) + FK)^T & W_2 \\
                   \end{array}
                 \right]
\succeq0.
%&\quad\quad W_1 = W_1^T,\quad W_2 = W_2^T,\nonumber\\
%&\quad\quad K_{SI}=K_{SI}^T,\;\;(1-\varepsilon)I_p+K_{SI}\succeq0.
\end{align}
We have relaxed $\tilde{\mathbb{P}}_{0}$ into the SDP \eqref{relaxation 2}. There are several types of efficient algorithms for solving SDPs, e.g., interior point methods and bundle method \citep{LV-SB:96}. These methods are implemented in commercial solvers such as Mosek, SeDuMi, and CVX, and can output the value of the SDP up to an additive error $\epsilon$ in time that is polynomial in the program description size and $\log{1}/{\epsilon}$.

\section{Problem $\tilde{\mathbb{P}}_2$ \label{Section P2}}

In the last section, we provide an SDP relaxation for problem $\tilde{\mathbb{P}}_{0}$.
%in which the $\ell_1$ norm heuristic, the nuclear norm heuristic and a positive semidefinite condition are applied to relax the $\ell_0$ norm minimization, the rank constraint and the determinant constraint, respectively.
This approach can be applied to problem $\tilde{\mathbb{P}}_2$ by replacing the $\ell_1$ norm heuristic by the SDP characterization of $\ell_2$ norm (page 170 of \citep{Boyd.Vandenberghe:04}) in \eqref{relaxation 2} and the resulting problem is:
\begin{align}
\label{P2 method one}
&\mathop {\min }\limits_{\begin{array}{*{20}c}
   \substack{ z\in\mathbb{C},K\in\mathbb{R}^{(p+l)\times(n+p)},t\in\mathbb{R} \\ K_{SI}\in\mathbb{S}^{p}, W_1  \in \mathbb{S}^{n + l},  W_2  \in \mathbb{S}^{n + p} }   \\
 \end{array} }t+c({\rm Tr}(W_1)+{\rm Tr}(W_2))\nonumber\\
&{\rm s.t.}\;  \left[ {\begin{array}{*{20}c}
   {tI_{p+l}} & [\bar H,\Pi]K  \\
   {K^T[\bar H,\Pi]^T } & {tI_{n+p}}  \\
 \end{array} } \right] \succeq 0,\;(1-\varepsilon)I_p+K_{SI}\succeq0\nonumber\\
&\quad\;\;\left[\begin{array}{cc}
                     W_1 & \bar D(z) + FK \\
                     (\bar D(z) + FK)^T & W_2 \\
                   \end{array}
                 \right]
\succeq0.
%&\quad\quad W_1 = W_1^T,\quad W_2 = W_2^T,\nonumber\\
%&\quad\quad K_{SI}=K_{SI}^T,\;\;(1-\varepsilon)I_p+K_{SI}\succeq0.
\end{align}
%This approach requires to tune the parameter $c$. If $c$ is too small, the perturbed system may still be strongly observable; if $c$ is too large, the matrix pencil in \eqref{the3} may lose rank much larger than one, which implies that the system is perturbed much more than necessary and the resulted $\|K\|_{0,\mathcal{L}}$ is likely to be large.
%This approach requires to tune the parameter $c$. The optimal $c$ can only be obtained experimentally. In this section, by assuming that the allowed perturbations restricted by $\mathbb{K}$ will always remain controllability, we investigate the scenario of $\mathbb{P}_2$ introduced in Section \ref{OP formulation} by equivalently formulating the constraint of $\mathbb{K}$ in a linearly constrained manner and provide a different approach to obtain the perturbation $K$ which does not need parameter tuning and can guarantee the non-strong-observability constraint.
%A limitation of this approach is that the rank deficiency constraint ${\rm rank}(\bar D(z)+FK)<\rho$ is relaxed to a soft constraint by the nuclear norm heuristic and the parameter $c$ needs to be empirically tuned.
For problem \eqref{P2 method one}, $c$ can only be tuned empirically. It is challenging to estimate the total time one needs to tune $c$ \emph{a priori}. For each given $c$, one needs to numerically solve the SDP of problem \eqref{P2 method one}. In this section, we study an approach which can analytically construct a feasible perturbation matrix $K$ that satisfies the constraint ${\rm rank}(\bar D(z)+FK)<\rho$ for a subclass of $\tilde{\mathbb{P}}_2$. This approach is more systematic as one can determine the largest possible tuning time of $\rho$ \emph{a priori}. Moreover, this approach is computationally more efficient than numerically solving the SDP of problem \eqref{P2 method one}.
%Notice that $\bar D(z)\in\mathbb{C}^{(n+q)\times(n+p)}$ and thus ${\rm rank}(\bar D(z)+FK)\leq\min\{n+q,n+p\}$. Hence, one can start with the maximum number $\rho=\min\{n+q,n+p\}$. If some target entries are not protected, one can then decrease $\rho$ by one, until all the target entries are protected or $\rho=1$. Hence, $\rho$ needs to be tuned for at most $\min\{n+q,n+p\}$ times. Moreover, for each given $\rho$, the approach only needs to perform a singular value decomposition (SVD) operation, which is computationally much cheaper than numerically solving the SDP of problem \eqref{P2 method one}.
In particular, in this section, we consider the following subclass of problem \eqref{relaxed OP agent-wise} where the controllability constraint is dropped:
\begin{align}
&\min\nolimits_{K\in\mathbb{R}^{(p+l)\times (n+p)},z\in\mathbb{C}}\left\|[\bar H,\Pi]K\right\|_2\nonumber\\
&{\rm s.t.}\;\; {\rm rank}(\bar D(z)+FK) < \rho.
\label{e6}
\end{align}
%where $\mathbb{K}=\mathbb{R}^{(p+l)\times(n+p)}$ and the controllability constraint is dropped.

\begin{remark}
\label{remark P2 controllability}
We next identify a class of problems where the perturbations do not affect system controllability so that problem \eqref{e6} can be applied. We rewrite system \eqref{state eq} and \eqref{e9} in the following form: $x(k+1)=\bar Ax(k)+\bar B^eu^e(k)+\bar B^cu^c(k)$ and $y(k)=\bar Gx(k)+\bar H^eu^e(k)+\bar H^cu^c(k)$, where $u^c$ is the control input while $u^e$ is some exogenous signal which is not used to control the system. Hence, we only need the system to be controllable with respect to $u^c$, rather than $u^e$. Assume that the target entries only include the entries of $u^e$ but do not include any entry of $u^c$. In this case, to protect privacy, we only need to perturb $(\bar A,\bar B^e,\bar G,\bar H^e)$, but do not need to perturb $(\bar B^c,\bar H^c)$. Assume that $\bar B^c$ has full row rank. Then, for any perturbed matrix $\hat A$, the controllability matrix with respect to $u^c$, $[\bar B^c,\hat A\bar B^c,\cdots,\hat A^{n-1}\bar B^c]$, always has full row rank, which implies that the perturbed system is always controllable with respect to $u^c$. Hence, for the above scenario, the perturbations do not affect system controllability with respect to the control inputs and problem \eqref{e6} can be applied. An example of the above scenario is the heating, ventilation, and air conditioning (HVAC) system in Section \ref{simulation section}.\oprocend
%In this example, the numbers of occupants are private exogenous signals and the temperatures of supplied air are non-private control inputs.\oprocend
%The above scenario holds for a heating, ventilation, and air conditioning (HVAC) system, in which the target entries are occupant numbers, which are not control inputs of the system. Please refer to Section \ref{section HVAC privacy} for the detailed description of the privacy issue of HVAC systems.
\end{remark}

%Problem \eqref{e6} can be applied to the cases where one does not have to maintain controllability with respect to the target entries in the perturbed system. This is the case for a heating, ventilation, and air conditioning (HVAC) system, in which the target entries are numbers of occupants of building zones, which are not control inputs of the system. Please refer to Section \ref{section HVAC privacy} for the detailed description of the privacy issue of HVAC systems.

In this section, we impose the following assumption.

\begin{assumption}
$\bar D(z)$ has full row rank $\forall z\in\mathbb{C}$.\oprocend
\label{asm: q equals p}
\end{assumption}

\begin{remark}
\label{remark full row rank}
%Notice that $\bar D(z)\in\mathbb{C}^{(n+q)\times(n+p)}$ and hence
Assumption \ref{asm: q equals p} implies $q\leq p$ and $\bar D(z)\bar D(z)^\dag=I_{n+q}$. Assumption \ref{asm: q equals p} can be efficiently checked as follows. Let
%$r_{\bar D}$ be the normal rank of $\bar D$, i.e.,
$r_{\bar D}=\max_{z\in\mathbb{C}}{\rm rank}(\bar D(z))$. For any $z'\in\mathbb{C}$, if ${\rm rank}(\bar D(z'))<r_{\bar D}$, then $z'$ is called an invariant zero of $\bar D$. Assumption \ref{asm: q equals p} is equivalent to that $r_{\bar D}=n+q$ and $\bar D$ does not have an invariant zero. Given $(\bar A,\bar B,\bar G,\bar H)$, to check whether Assumption \ref{asm: q equals p} holds, one can first check whether $\bar D$ has an invariant zero. There are efficient algorithms to compute invariant zeros \citep{AEN-PVD:1982}
%\citep{AEN-PVD:1982,PM-PVD-AV:1994}
~\footnote{The algorithm of \citep{AEN-PVD:1982} for finding the invariant zeros of a linear system has been implemented in Matlab by the command $tzero$.}. If $\bar D$ does not have an invariant zero, then one can derive the value of $r_{\bar D}$ by computing $r_{\bar D}={\rm rank}(\bar D(z))$ with any $z\in\mathbb{C}$ and then check whether $r_{\bar D}=n+q$.

A sufficient condition for Assumption \ref{asm: q equals p} is that $q=p$ and system $(\bar A,\bar B,\bar G,\bar H)$ is strongly observable. By Theorem 3 and Corollary 4 of \citep{SY-MZ-EF:2016}, $(\bar A,\bar B,\bar G,\bar H)$ is strongly observable if and only if ${\rm rank}(\bar D(z))=n+p$ for all $z\in\mathbb{C}$. Hence, if $q=p$ and $(\bar A,\bar B,\bar G,\bar H)$ is strongly observable, we have ${\rm rank}(\bar D(z))=n+p=n+q$ for all $z\in\mathbb{C}$ and thus Assumption \ref{asm: q equals p} holds.\hfill\oprocend
\end{remark}

\subsection{Feasible solution for fixed $z$}

To solve problem \eqref{e6}, we first fix any $z\in \mathbb{C}$ and consider the following problem $\tilde{\mathbb{P}}_2(z)$:
\begin{align}
\label{fixed z P2}
&\min\nolimits_{K\in\mathbb{R}^{(p+l)\times (n+p)}}\left\|[\bar H,\Pi]K\right\|_2\nonumber\\
&{\rm s.t.}\;\; {\rm rank}(\bar D(z)+FK) < \rho.
\end{align}
We aim to find a feasible solution for problem \eqref{fixed z P2} and derive an upper bound of the optimal value of problem \eqref{fixed z P2}. To do this, we perform the singular value decomposition (SVD) on $\bar D(z)^\dag F$ as $\bar D(z)^\dag F=U(z)\Sigma(z) V(z)^T$. Since $\bar D(z)$ is assumed to have full row rank, we have ${\rm rank}(\bar D(z)^\dag F)={\rm rank}(F)$. For any integer $\ell$ such that $1\leq\ell\leq \min\{n+p,p+l\}$, let $\sigma_\ell(z)$ be the $\ell$-th diagonal entry of $\Sigma(z)$. Without loss of generality, assume that the $\sigma_\ell(z)$'s are arranged in the descending order, i.e., $\sigma_1(z)\geq\sigma_2(z)\geq\cdots\geq\sigma_{{\rm rank}(F)}(z)>0$ and $\sigma_{{\rm rank}(F)+1}(z)=\sigma_{{\rm rank}(F)+2}(z)=\cdots=\sigma_{\min\{n+p,p+l\}}(z)=0$. For any integer $\ell$ such that $1\leq\ell\leq n+p$, let $u_\ell(z)$ be the $\ell$-th column vector of $U(z)$. For any integer $\ell$ such that $1\leq\ell\leq p+l$, let $v_\ell(z)$ be the $\ell$-th column vector of $V(z)$. We then have $\bar D(z)^\dag F=\sum_{\ell=1}^{{\rm rank}(F)}\sigma_\ell(z) u_\ell(z) v_\ell(z)^T$.
%We have the following relations:
%\begin{align}
%\label{eq00_00}
%&F^\dag D(z)=\sum_{\ell=1}^{r(z)}\sigma_\ell(z)u_\ell(z)v_\ell(z)^T,\\
%%\label{eq00_1}
%%&v_\ell(z)^Tv_\ell(z)=1,\;\forall 1\leq\ell\leq n+p\\
%%&u_\ell(z)^TD(z)^{-1} Fv_\ell(z)=\sigma_\ell(z),\\
%\label{eq00_2}
%&v_\ell(z)^Tv_{\ell'}(z)=0,\;\forall 1\leq\ell,\ell'\leq n+p\;{\rm and}\;\ell\neq\ell'.
%%&u_\ell(z)^TD(z)^{-1} Fv_{\ell'}(z)=0,\;\forall\ell\neq\ell'.
%\end{align}

\begin{lemma}
\label{lemma P2 upper bound}
Suppose that Assumption \ref{asm: q equals p} holds. Fix any $z\in\mathbb{C}$. The following statements hold:\\
(i) Problem \eqref{fixed z P2} is feasible if and only if $\rho\geq n+q-{\rm rank}(F)+1$.\\
(ii) If $\rho>n+q$, then $\tilde K=\textbf{0}_{(p+l)\times(n+p)}$ is an optimal solution of problem \eqref{fixed z P2} and the optimal value of problem \eqref{fixed z P2} is zero.\\
(iii) If $n+q-{\rm rank}(F)+1\leq\rho\leq n+q$, then $\tilde K(z)=-\sum_{\ell=1}^{n+q-\rho+1}\frac{1}{\sigma_\ell(z)}v_\ell(z)u_\ell(z)^T$ is a feasible solution of problem \eqref{fixed z P2} and the optimal value of problem \eqref{fixed z P2} is upper bounded by $\|[\bar H,\Pi]\|_2\|(\bar D(z)^\dag F)^\dag\|_*$.\hfill\oprocend
\end{lemma}

\textbf{Proof:} We first prove (i). Given any matrix $M$, denote by $\mathcal{N}(M)$ the right null space of $M$ and denote by $|\mathcal{N}(M)|$ the dimension of $\mathcal{N}(M)$. First, we show that problem \eqref{fixed z P2} is infeasible if $\rho<n+q-{\rm rank}(F)+1$. We show this by contradiction. Given any $\rho<n+q-{\rm rank}(F)+1$. Assume that there exists some $K$ such that ${\rm rank}(\bar D(z)+FK)<\rho$. This implies
\begin{align}
\label{number of right null}
|\mathcal{N}(\bar D(z)+FK)|%\geq n+p-\rho+1
>p-q+{\rm rank}(F).
\end{align}
Let $u$ be any right null vector of $\bar D(z)+FK$, i.e., $(\bar D(z)+FK)u=\textbf{0}_{n+q}$. We then have $\bar D(z)u=-FKu$. There can only be two cases: (a) $\bar D(z)u=-FKu=\textbf{0}_{n+q}$ and (b) $\bar D(z)u=-FKu\neq\textbf{0}_{n+q}$. For case (a), $\bar D(z)u=\textbf{0}_{n+q}$ implies that $u\in\mathcal{N}(\bar D(z))$. By Assumption \ref{asm: q equals p}, ${\rm rank}(\bar D(z))=n+q$. Hence, $|\mathcal{N}(\bar D(z))|=n+p-{\rm rank}(\bar D(z))=n+p-(n+q)=p-q$. Thus, the dimension of the space of $u$ for case (a) is $p-q$. For case (b), $-FKu\neq\textbf{0}_{n+q}$ implies that $u$ is in the complementary space of $\mathcal{N}(\bar D(z)+FK)$ and thus the dimension of the space of $u$ in this case equals to ${\rm rank}(FK)\leq{\rm rank}(F)$. Combining the two cases (a) and (b), we reach that $|\mathcal{N}(\bar D(z)+FK)|\leq p-q+{\rm rank}(F)$, which contradicts \eqref{number of right null}. Hence, problem \eqref{fixed z P2} is infeasible if $\rho<n+q-{\rm rank}(F)+1$.

We next show that problem \eqref{fixed z P2} is feasible if $\rho\geq n+q-{\rm rank}(F)+1$. This is proven by proving (ii) and (iii). We first prove (ii). Since ${\rm rank}(\bar D(z))=n+q$, ${\rm rank}(\bar D(z))<\rho$. It is clear that $\tilde K=\textbf{0}_{(p+l)\times(n+p)}$ is a feasible solution for problem \eqref{fixed z P2}. Since $\|[\bar H,\Pi]K\|_2\geq0$ and $\|[\bar H,\Pi]\tilde K\|_2=0$ with $\tilde K=\textbf{0}_{(p+l)\times(n+p)}$, we have that $\tilde K=\textbf{0}_{(p+l)\times(n+p)}$ is an optimal solution for problem \eqref{fixed z P2} and the optimal value of problem \eqref{fixed z P2} is zero.

We next prove (iii). Substituting $\tilde K(z)=-\sum_{\ell=1}^{n+q-\rho+1}$ $\frac{1}{\sigma_\ell(z)}v_\ell(z)u_\ell(z)^T$ into $\bar D(z)+FK$ yields:
\begin{align}
\label{perturbed pencil matrix}
&\bar D(z)+F\tilde K(z)=\bar D(z)(I_{n+p}+\bar D(z)^\dag F\tilde K(z))\nonumber\\
%&=\bar D(z)(I_{n+p}-\bar D(z)^\dag F\sum_{\ell=1}^{n+q-\rho+1}\frac{1}{\sigma_\ell(z)}v_\ell(z)u_\ell(z)^T)\nonumber\\
&=\bar D(z)(I_{n+p}-\sum_{\ell=1}^{n+q-\rho+1}\frac{1}{\sigma_\ell(z)}\bar D(z)^\dag Fv_\ell(z)u_\ell(z)^T)\nonumber\\
&=\bar D(z)(I_{n+p}-\sum\nolimits_{\ell=1}^{n+q-\rho+1}u_\ell(z)u_\ell(z)^T)\nonumber\\
&=\bar D(z)(\sum_{\ell=1}^{n+p}u_\ell(z)u_\ell(z)^T-\sum_{\ell=1}^{n+q-\rho+1}u_\ell(z)u_\ell(z)^T)\nonumber\\
&=\bar D(z)\sum\nolimits_{\ell=n+q-\rho+2}^{n+p}u_\ell(z)u_\ell(z)^T.
\end{align}
In \eqref{perturbed pencil matrix}, the first equation is a result of Assumption \ref{asm: q equals p}; the fourth equation is due to $\bar D(z)^\dag Fv_\ell(z)=\sigma_\ell(z)u_\ell(z)$ for all $\ell\in[1,{\rm rank}(F)]$; and the fifth equation is because $U(z)$ is unitary, i.e., $U(z)U(z)^T=\sum_{\ell=1}^{n+p}u_\ell(z)u_\ell(z)^T=I_{n+p}$. Since $U(z)$ is unitary, $u_\ell(z)^Tu_{\ell'}(z)=0$ for all $\ell,\ell'\in[1,n+p]$ and $\ell\neq\ell'$. By \eqref{perturbed pencil matrix}, for all $\ell'\in[1,n+q-\rho+1]$, we have\\
$(\bar D(z)+F\tilde K(z))u_{\ell'}(z)\\
=\bar D(z)\sum\nolimits_{\ell=n+q-\rho+2}^{n+p}u_\ell(z)u_\ell(z)^Tu_{\ell'}(z)=\textbf{0}_{n+q}$.\\
As ${\rm rank}(\bar D(z))=n+q\leq n+p$ and $U(z)$ is unitary, there must be exactly $p-q$ column vectors of $U(z)$ which are right null vectors of $\bar D(z)$. Let $u_{\ell'}(z)$ be any such vector, i.e., $u_{\ell'}(z)$ is a column vector of $U(z)$ and $\bar D(z)u_{\ell'}(z)=\textbf{0}_{n+q}$. We show that $\ell'\geq{\rm rank}(F)+1$, i.e., $u_{\ell'}(z)$ corresponds to a zero singular value of $\bar D(z)^\dag F$. As $\bar D(z)u_{\ell'}(z)=\textbf{0}_{n+q}$, $u_{\ell'}(z)^T\bar D(z)^T=\textbf{0}_{1\times(n+q)}$. So $u_{\ell'}(z)^T\bar D(z)^\dag=\textbf{0}_{1\times(n+q)}$. Hence, $u_{\ell'}(z)^T\bar D(z)^\dag F=\sigma_{\ell'}(z)v_{\ell'}(z)^T=\textbf{0}_{1\times(p+l)}$. Thus $\sigma_{\ell'}(z)=0$ and $\ell'\in[{\rm rank}(F)+1,n+p]$. As $\rho\geq n+q-{\rm rank}(F)+1$, ${\rm rank}(F)+1\geq n+q-\rho+2$. By \eqref{perturbed pencil matrix}, we have\\
$(\bar D(z)+F\tilde K(z))u_{\ell'}(z)=\bar D(z)\sum\nolimits_{\ell=n+q-\rho+2}^{n+p}u_\ell(z)\\ \cdot u_\ell(z)^Tu_{\ell'}(z)=\bar D(z)u_{\ell'}(z)=\textbf{0}_{n+q}$.\\
We have found $(n+q-\rho+1)+(p-q)=n+p-\rho+1$ linearly independent right null vectors of $\bar D(z)+F\tilde K$. Hence, ${\rm rank}(\bar D(z)+F\tilde K)<\rho$ and $\tilde K(z)$ $=-\sum_{\ell=1}^{n+q-\rho+1}\frac{1}{\sigma_\ell(z)}v_\ell(z)u_\ell(z)^T$ is a feasible solution of problem \eqref{fixed z P2}. Since $\rho\geq n+q-{\rm rank}(F)+1$, we have ${\rm rank}(F)\geq n+q-\rho+1$. Hence, for the constructed $\tilde K(z)$, the value of the objective function of problem \eqref{fixed z P2} is
\begin{align}
\label{P2i objective function 2norm}
&\!\!\!\!\|[\bar H,\Pi]\tilde K(z)\|_2=\|\sum\nolimits_{\ell=1}^{n+q-\rho+1}\frac{[\bar H,\Pi]}{\sigma_\ell(z)}v_\ell(z)u_\ell(z)^T\|_2\nonumber\\
%=\frac{\|\Psi_\kappa v_\kappa(z)u_\kappa(z)^T\Phi_\kappa\|_2}{\|\Phi_\kappa D(z)^{-1} F\Psi_\kappa\|_2}\nonumber\\
&\!\!\!\!\leq\sum\nolimits_{\ell=1}^{{\rm rank}(F)}\frac{\|[\bar H,\Pi]\|_2}{\sigma_\ell(z)}{\|v_\ell(z)u_\ell(z)^T\|_2}\nonumber\\
%&\leq\sum_{\ell=1}^{n+p-\rho+1}\frac{\|v_\ell(z)\|_2\|u_\ell(z)\|_2}{\sigma_\ell(z)}\nonumber\\
&\!\!\!\!=\sum\nolimits_{\ell=1}^{{\rm rank}(F)}\frac{\|[\bar H,\Pi]\|_2}{\sigma_\ell(z)}=\|[\bar H,\Pi]\|_2\|(D(z)^\dag F)^\dag\|_*.
%&\leq\frac{\|\Psi_\kappa\|_2\|v_\kappa(z)u_\kappa(z)^T\|_2\|\Phi_\kappa\|_2}{\|\Phi_\kappa D(z)^{-1} F\Psi_\kappa\|_2}\nonumber\\
%&=\|\Psi_\kappa\|_2\|\Phi_\kappa\|_2\|\Phi_\kappa D(z)^{-1} F\Psi_\kappa\|_2^{-1}\nonumber\\
%&=\|\Phi_\kappa D(z)^{-1} F\Psi_\kappa\|_2^{-1}.
\end{align}
Hence, $\|[\bar H,\Pi]\|_2\|(\bar D(z)^\dag F)^\dag\|_*$ is an upper bound of the optimal value of problem \eqref{fixed z P2}.\oprocend

\begin{remark}
%Notice that $\sigma_1(z)=\|D(z)^{-1}F\|$. Hence, in the case where $\rho=n+p$, we can see that the upper bound given by Lemma \ref{lemma P2 upper bound} reduces to the lower bound given by Lemma \ref{lemma P2 lower bound}, which implies that the perturbation matrix $K(z)$ constructed by \eqref{Ki choice} achieves optimality for problem \eqref{fixed z P2} in this special case.
%Since $D(z)$ is square and invertible, it holds that $D(z)D(z)^{-1}=I_{n+p}$ and thus $D(z)c_\kappa(z)=D(z)D(z)^{-1}F\Psi_\kappa v_\kappa(z)=F\Psi_\kappa v_\kappa(z)$, which cancels $D(z)$ and enables the rest steps of \eqref{eq use inverse property}.
%The first equation of \eqref{eq use inverse property} makes use of Assumption \ref{asm: q equals p}. Without Assumption \ref{asm: q equals p}, we need to use the pseudo inverse of $D(z)$ and $c_\ell(z)=D(z)^\dag Fv_\ell(z)$. However, since $D(z)$ has more rows than columns, $D(z)D(z)^\dag\neq I_{n+q}$ (instead, $D(z)^\dag D(z)=I_{n+p}$). Thus, $D(z)$ cannot be cancelled and the rest steps of \eqref{eq use inverse property} do not follow.
The derivation of \eqref{P2i objective function 2norm} utilizes the property $\|v_\ell(z)u_\ell(z)^T\|_2=1$. This property in general does not hold for $\ell_0$ norm or $\ell_1$ norm, i.e., $\|v_\ell(z)u_\ell(z)^T\|_0\neq1$ and $\|v_\ell(z)u_\ell(z)^T\|_1\neq1$. Thus, the approach developed in this section is not suitable for problem $\tilde{\mathbb{P}}_{0}$.\oprocend
\label{remark: optimal P2z}
\end{remark}

\subsection{Minimization over $z$\label{section z minimization}}

To compute the feasible perturbation matrix by Lemma \ref{lemma P2 upper bound}, we first need to decide the value of $z$.
%This problem is studied in this subsection.
By Lemma \ref{lemma P2 upper bound}, $\|[\bar H,\Pi]\|_2\|(\bar D(z)^\dag F)^\dag\|_*$ is an upper bound of the optimal value of problem \eqref{fixed z P2}. Since $\|[\bar H,\Pi]\|_2$ is a constant, we aim to minimize $\|(\bar D(z)^\dag F)^\dag\|_*$ over $z\in\mathbb{C}$. Notice that $\bar D(z)$ is a symbolic matrix and computing the pseudo inverse of a symbolic matrix could be time-consuming. To avoid it, we perform the following relaxation.

We perform SVD as $F=U_F\Sigma_FV_F$. Let $r={\rm rank}(F)$ and $\Sigma_F=\left[ {\begin{array}{*{20}{c}}
{{{\tilde \Sigma }_F}}&\textbf{0}_{r\times(p+l-r)}\\
\textbf{0}_{(n+q-r)\times r}&\textbf{0}_{(n+q-r)\times(p+l-r)}
\end{array}} \right]$, where $\tilde\Sigma_F$ is a diagonal matrix whose diagonal entries are the $r$ non-zero singular values of $F$. Let $U_F=[U_{F1},U_{F2}]$ and $V_F=[V_{F1},V_{F2}]$, with $U_{F1}\in\mathbb{R}^{(n+q)\times r}$, $U_{F2}\in\mathbb{R}^{(n+q)\times(n+q-r)}$, $V_{F1}\in\mathbb{R}^{(p+l)\times r}$ and $V_{F2}\in\mathbb{R}^{(p+l)\times(p+l-r)}$. The SVD of $F$ is then
$F=\left[ {\begin{array}{*{20}{c}}
U_{F1}&U_{F2}
\end{array}} \right]$ $\left[ {\begin{array}{*{20}{c}}
{{{\tilde \Sigma }_F}}&\textbf{0}_{r\times(p+l-r)}\\
\textbf{0}_{(n+q-r)\times r}&\textbf{0}_{(n+q-r)\times(p+l-r)}
\end{array}} \right]\left[ {\begin{array}{*{20}{c}}
V_{F1}^T\\
V_{F2}^T
\end{array}} \right]=U_{F1}\tilde\Sigma_FV_{F1}^T$. Since $\bar D(z)^\dag U_{F1}$ has full column rank and $\tilde\Sigma_FV_{F1}^T$ has full row rank, by Corollary 1.4.2 of \citep{SLC-CDM:2009}, $(\bar D(z)^\dag F)^\dag=(\bar D(z)^\dag U_{F1}\tilde\Sigma_FV_{F1}^T)^\dag=(\tilde\Sigma_FV_{F1}^T)^\dag(\bar D(z)^\dag U_{F1})^\dag$. By the definitions of matrix $\ell_2$ norm and nuclear norm \citep{RH-CJ:1985}, we have $\|M\|_*\leq r_M\|M\|_2$ for any matrix $M$ with rank $r_M$. It follows that
\begin{align}
\label{revised z minimization eq 1}
&\|(\bar D(z)^\dag F)^\dag\|_*\leq r\|(\bar D(z)^\dag F)^\dag\|_2\nonumber\\
&\leq r\|(\tilde\Sigma_FV_{F1}^T)^\dag\|_2\|(\bar D(z)^\dag U_{F1})^\dag\|_2.
\end{align}
Since $\bar D(z)^\dag$ has full column rank and $U_F$ is unitary, by Ex. 2 in page 80 of \citep{ICFI:2009}, $\sigma_{\min}(\bar D(z)^\dag U_{F1})\geq\sigma_{\min}(\bar D(z)^\dag U_F)$. Since $\|(\bar D(z)^\dag U_{F1})^\dag\|_2=\frac{1}{\sigma_{\min}(\bar D(z)^\dag U_{F1})}$ and $\|(\bar D(z)^\dag U_F)^\dag\|_2=\frac{1}{\sigma_{\min}(\bar D(z)^\dag U_{F})}$, we then have $\|(\bar D(z)^\dag U_{F1})^\dag\|_2\leq\|(\bar D(z)^\dag U_F)^\dag\|_2$. By \eqref{revised z minimization eq 1}, we have\\
$\|(\bar D(z)^\dag F)^\dag\|_*\leq r\|(\tilde\Sigma_FV_{F1}^T)^\dag\|_2\|(\bar D(z)^\dag U_{F})^\dag\|_2\\
=r\|(\tilde\Sigma_FV_{F1}^T)^\dag\|_2\|U_F^T\bar D(z)\|_2=r\|(\tilde\Sigma_FV_{F1}^T)^\dag\|_2\|\bar D(z)\|_2\\
\leq r\|(\tilde\Sigma_FV_{F1}^T)^\dag\|_2\|\bar D(z)\|_F$.\\
The first equality is because $\bar D(z)^\dag$ has full column rank and $U_F$ has full row rank (see Corollary 1.4.2 in page 22 of \citep{SLC-CDM:2009}); the second equality holds because $U_F$ is a unitary matrix; and the last inequality is due to the equivalence of matrix norms \citep{RH-CJ:1985}, i.e., $\|M\|_2\leq\|M\|_F$ for any matrix $M$. Notice that $r\|(\tilde\Sigma_FV_{F1}^T)^\dag\|_2$ is independent of $z$. Hence, we are to minimize $\|\bar D(z)\|_F$ or equivalently $\|\bar D(z)\|_F^2$ over $z\in\mathbb{C}$.
%Notice that $\bar D(z)^T\bar D(z)=\left[ {\begin{array}{*{20}{c}}
%z^2I_{n}-(\bar A+\bar A^T)z+\bar A^T\bar A+\bar G^T\bar G&\bar G^T\bar H-(zI_n-\bar A^T)\bar B\\
%\bar H^T\bar G-\bar B^T(zI_n-\bar A)&\bar B^T\bar B+\bar H^T\bar H
%\end{array}} \right]$.
We have\\
$\|\bar D(z)\|_F^2={\rm Tr}(\bar D(z)^T\bar D(z))={\rm Tr}(z^2I_{n}-(\bar A+\bar A^T)z)\\
+{\rm Tr}(\bar A^T\bar A+\bar G^T\bar G)+{\rm Tr}(\bar B^T\bar B+\bar H^T\bar H)$.\\
Since ${\rm Tr}(\bar A^T\bar A+\bar G^T\bar G)+{\rm Tr}(\bar B^T\bar B+\bar H^T\bar H)$ is constant, we are to minimize ${\rm Tr}(z^2I_{n}-(\bar A+\bar A^T)z)=nz^2-2{\rm Tr}(\bar A)z$. Hence, the optimal value of $z$ is $\tilde z={{\rm Tr}(\bar A)}/{n}$.

\subsection{Overall approach}

Given $\rho\geq n+q-{\rm rank}(F)+1$, we have derived a procedure to determine a feasible solution of problem \eqref{e6} which minimizes an upper bound of the optimal value of \eqref{e6}. The procedure is summarized in Algorithm \ref{Algorithm P2}.

\begin{algorithm}[htbp]
\small
\caption{Algorithm for finding a suboptimal feasible solution of problem \eqref{e6}} \label{Algorithm P2}

Compute $\tilde z={{\rm Tr}(\bar A)}/{n}$;

Perform SVD: $\bar D(\tilde z)^\dag F=\sum_{\ell=1}^{{\rm rank}(F)}\sigma_\ell(\tilde z) u_\ell(\tilde z) v_\ell(\tilde z)^T$;

Compute $\tilde K$ by (ii) or (iii) of Lemma \ref{lemma P2 upper bound}:\\
%\nonl
if $\rho>n+q$, $\tilde K=\textbf{0}_{(p+l)\times(n+p)}$;\\
%\nonl
if $\rho\leq n+q$, $\tilde K=-\sum_{\ell=1}^{n+q-\rho+1}\frac{1}{\sigma_\ell(\tilde z)}v_\ell(\tilde z)u_\ell(\tilde z)^T$.
\end{algorithm}
\normalsize

As mentioned in Section \ref{rank relaxation section}, we aim to protect all the target entries with the largest possible $\rho$. The tuning of $\rho$ can be systematically performed as follows. By Assumption \ref{asm: q equals p}, we have ${\rm rank}(\bar D(z)+FK)\leq\min\{n+q,n+p\}=n+q$. Hence, we can start with the maximum number $\rho=n+q$ and run Algorithm \ref{Algorithm P2}. After that, we use the mechanism introduced at the second last paragraph of Section \ref{rank relaxation section} to check whether all the target entries are protected. If not, we decrease $\rho$ by one and re-run Algorithm \ref{Algorithm P2}. The procedure is repeated until all the target entries are protected or $\rho=1$. Hence, $\rho$ needs to be tuned for at most $n+q$ times. For each given $\rho$, most computational effort of Algorithm \ref{Algorithm P2} is spent to perform the SVD operation, which is computationally more efficient than numerically solving the SDP of problem \eqref{P2 method one}.

\section{Case study\label{simulation section}}

In this section, we validate the efficacy of the developed techniques by an HVAC system.

\subsection{System model}

\begin{table}[!t]
\renewcommand{\arraystretch}{1}
\caption{Parameters/variables of the HVAC system}
\label{thermal parameters}
\scriptsize
\centering
\begin{tabular}{|l|l|}
\hline
$L_i$ & thermal capacity of zone $\mathcal{Z}_i$\\
\hline
$R_{ji}$ & thermal conductance between zone $\mathcal{Z}_i$ and zone $\mathcal{Z}_j$\\
\hline
$\Delta t$ & discretization stepsize\\
\hline
$c_o$ & thermal load per occupant\\
\hline
$T_i^s$ & temperature of air supplied to zone $\mathcal{Z}_i$\\
\hline
$T_i$ & temperature of zone $\mathcal{Z}_i$\\
\hline
$m_i^s$ & mass flow rate of air supplied to zone $\mathcal{Z}_i$\\
\hline
$c_p$ & thermal capacity of air\\
\hline
$V_i$ & number of occupants of zone $\mathcal{Z}_i$\\
\hline
\end{tabular}
\end{table}
\normalsize

Consider a set of $N$ building zones $\mathcal{Z}=\{\mathcal{Z}_1,\cdots,\mathcal{Z}_N\}$. The physical meanings of the system parameters and variables are listed in Table \ref{thermal parameters}. For each $i\in\{1,\cdots,N\}$, the following discrete-time dynamic model of zone $\mathcal{Z}_i$ is adopted from \citep{AK-FB:2011}:
%\begin{align}
%\label{HVAC}
%&L_i\frac{T_i(k+1)-T_i(k)}{\Delta t}\nonumber\\
%&=\sum\limits_{j\in\mathcal{N}_i}R_{ji} \left(T_j (k)-\frac{T_i(k+1)+T_i(k)}{2}\right)+c_oV_i(k)\nonumber\\
%&+m_i^s(k)c_p\left(T_i^s(k)-\frac{T_i(k+1)+T_i(k)}{2}\right).
%\end{align}
\begin{align}
\label{HVAC}
&({L_i}/{\Delta t}+\sum\nolimits_{j\in\mathcal{N}_i}R_{ji}/2+m_i^s(k)c_p/2)T_i(k+1)\nonumber\\
&=({L_i}/{\Delta t}-\sum\nolimits_{j\in\mathcal{N}_i}R_{ji}/2-m_i^s(k)c_p/2)T_i(k)\nonumber\\
&+\sum\nolimits_{j\in\mathcal{N}_i}R_{ji}T_j(k)+m_i^s(k)c_pT_i^s(k)+c_oV_i(k).
\end{align}
Assume that $m_i^s$ is constant, i.e., $m_i^s(k)\equiv \bar m_i^s$ for all $i$'s. The state and the control input of each zone $\mathcal{Z}_i$ is $T_i$ and $T_i^s$, respectively. For each $i\in\{1,\cdots,N\}$, let $x_i(k)=T_i(k)$, $u_i^e(k)=V_i(k)$ and $u_i^c(k)=T_i^s(k)$. System \eqref{HVAC} can then be written as $x_i(k+1)=\bar A_{ii}x_i(k)+\sum\nolimits_{j\in\mathcal{N}_i}\bar A_{ij}x_j(k)+\bar B_i^eu_i^e(k)+\bar B_i^cu_i^c(k)$, where $\bar A_{ii}=\frac{L_i/\Delta t-\bar m_i^sc_p/2-1/2\sum_{j\in\mathcal{N}_i}R_{ji}}{L_i/\Delta t+\bar m_i^sc_p/2+1/2\sum_{j\in\mathcal{N}_i}R_{ji}}$,\\
$\bar A_{ij}=\frac{R_{ji}}{L_i/\Delta t+\bar m_i^sc_p/2+1/2\sum_{j\in\mathcal{N}_i}R_{ji}},\;\forall j\in\mathcal{N}_i$,\\
$\bar B_i^e=\frac{c_o}{L_i/\Delta t+\bar m_i^sc_p/2+1/2\sum_{j\in\mathcal{N}_i}R_{ji}}$,\\
$\bar B_i^c=\frac{m_i^sc_p}{L_i/\Delta t+\bar m_i^sc_p/2+1/2\sum_{j\in\mathcal{N}_i}R_{ji}}$. The outputs required by the data requester are given by $y(k)=\bar G x(k)+\bar H^e u^e(k)+\bar H^cu^c(k)$ and the data requester knows $(\bar A,\bar B^e,\bar B^c,\bar G,\bar H^e,\bar H^c)$. The above state and output equations can be written in the form of \eqref{state eq} and \eqref{e9} with $u=[u^{eT},u^{cT}]^T$, $\bar B=[\bar B^e,\bar B^c]$ and $\bar H=[\bar H^e,\bar H^c]$.

\subsection{Privacy issue\label{section HVAC privacy}}

The usage of occupancy data 
%\citep{BB-JX-AN-RG-YA:2013,VLE-AEC:2010,WK-SS-FM:2014}, 
poses risks on the privacy of individual occupants. It has been shown in \citep{XW-PT:2014} that with some auxiliary information such as an office directory, individual location traces can be inferred from the occupancy data with accuracy of more than $90\%$. The information attached to location traces could reveal much about the individual occupants' interests, activities and relationships \citep{ML-DM-SW:2010}.

In system \eqref{HVAC}, the individual location trace is the private information. As mentioned above, this information could potentially be inferred from the occupancy data $V_i$'s. We aim to use the proposed intentional input-output perturbations to perturb system \eqref{HVAC} so that the perturbed system is private in the sense of Definition \ref{def: agent-wise privacy}.
%and thus the data requester cannot recover the occupancy data from the historical system outputs.%Meanwhile, we require that the perturbed system remains controllable.

\subsection{Applicability of the developed techniques\label{section applicabiliity}}

Problem $\tilde{\mathbb{P}}_{0}$: In the above HVAC system, $u^c$ is the control while $u^e$ is an exogenous signal which is not used to control the system. Hence, when we formulate problem \eqref{relaxed OP agent-wise}, we should only maintain controllability with respect to $u^c$, but not $u^e$. This is embedded into problem \eqref{relaxed OP agent-wise} by replacing the input matrix $\bar B$ in the controllability constraint by the partial input matrix $\bar B^c$ associated with $u^c$. One can then apply the relaxation techniques proposed in Section \ref{P0 section}. The simulation results for problem $\tilde{\mathbb{P}}_{0}$ in this section are derived for the modified problem.

Problem $\tilde{\mathbb{P}}_{2}$: Notice that the target entries only include entries of $u^e$, but no entry of $u^c$. Moreover, in our problem, $\bar B^c$ has full row rank (the parameters of $\bar B^c$ are adopted from \citep{YM-GA-FB:2011}). As mentioned in Remark \ref{remark P2 controllability}, to protect privacy, we only need to perturb $(\bar A,\bar B^e,\bar G,\bar H^e)$, but do not need to perturb $(\bar B^c,\bar H^c)$, and it is guaranteed that the perturbed system is controllable with respect to $u^c$. Hence, problem \eqref{e6} can be applied. In the following simulation for $\tilde{\mathbb{P}}_{2}$, the matrices $\bar D(z)$ and $F$ in problem \eqref{e6} are defined by $(\bar A,\bar B^e,\bar G,\bar H^e)$. For the simulation for problem $\tilde{\mathbb{P}}_{0}$, in order to verify the relaxation for the controllability constraint (the constraint of problem \eqref{relaxation 1}) proposed in Section \ref{P0 section}, we perturb the overall matrices $(\bar A,\bar B,\bar G,\bar H)$ where $\bar B=[\bar B^e,\bar B^c]$ and $\bar H=[\bar H^e,\bar H^c]$. Accordingly, the matrices $\bar D(z)$ and $F$ in problem \eqref{relaxation 2} are defined by $(\bar A,\bar B,\bar G,\bar H)$.

\subsection{Simulation results\label{simulation results}}

\begin{figure}
\begin{center}
\includegraphics[width=1\linewidth]{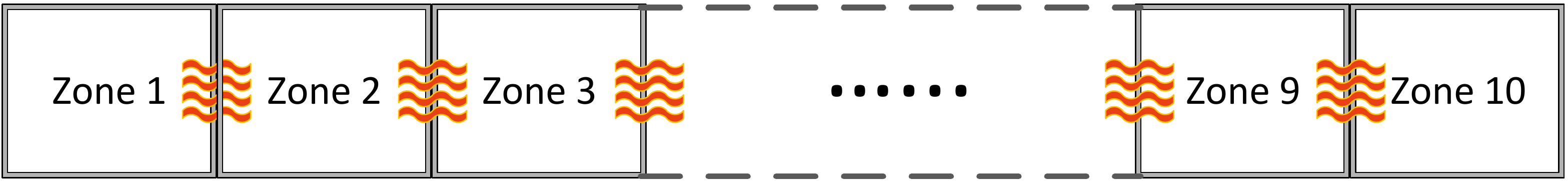}
\caption{Floor plan for the HVAC system}
\label{floor_plan}
\end{center}
\end{figure}

We take $N=10$, which leads to $n=10$ and $p^e=p^c=10$, where $p^e$ and $p^c$ are the dimensions of $u^e$ and $u^c$, respectively. We choose $q^e=q^c=10$, where $q^e$ and $q^c$ are the row numbers of $\bar H^e$ and $\bar H^c$, respectively. The undirected graph describing the topology of the zone network is denoted by $\mathcal{G}=(\mathcal{V},\mathcal{E})$, where $\mathcal{V}=\{\mathcal{Z}_1,\cdots,\mathcal{Z}_{10}\}$ and $\mathcal{E}=\{(\mathcal{Z}_1,\mathcal{Z}_2),\cdots,(\mathcal{Z}_9,\mathcal{Z}_{10})\}$. The floor plan is depicted by Fig. \ref{floor_plan}. This adjacency topology is widely used in the literature, e.g., \citep{YM-GA-FB:2011}. The values of the parameters in Table \ref{thermal parameters} are adopted from \citep{YM-GA-FB:2011}.
%such that $q=p$ and the unperturbed system is strongly observable and thus Assumption \ref{asm: q equals p} holds (see Remark \ref{remark full row rank}).

%\subsubsection{Problem $\tilde{\mathbb{P}}_{0,\mathcal{L}}$}

\begin{figure}[H]
\begin{center}
\includegraphics[width=1\linewidth]{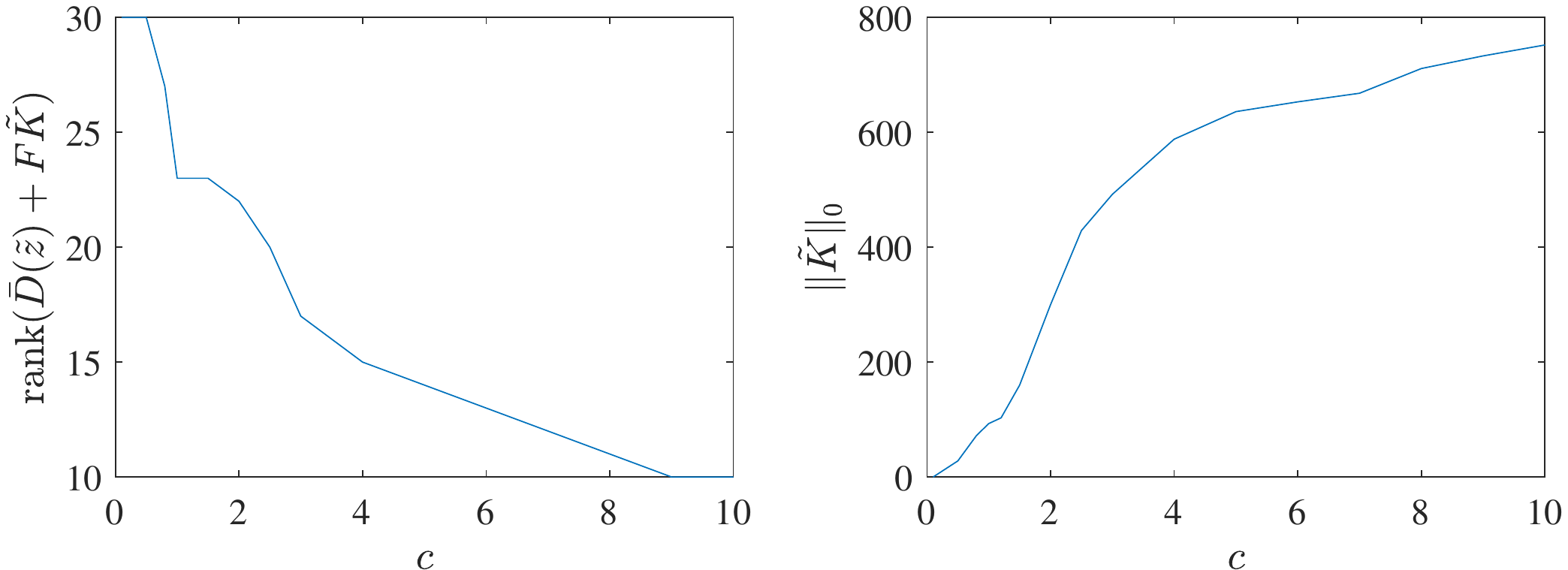}
\caption{Effect of $c$ with $\varepsilon\equiv0.1$ for problem $\tilde{\mathbb{P}}_{0}$}
\label{P0 figure}
\end{center}
\end{figure}

\textbf{Problem $\tilde{\mathbb{P}}_{0}$}. 
%we construct the constraint set $\mathbb{K}$ by randomly choosing $100$ positions where we cannot install sensors or communication links. We have examined that, under the structure of $\mathbb{K}$, there exist realizations of $K$ such that the perturbed system is uncontrollable so that the controllability constraint is nontrivial. 
The matrices $\bar G$ and $\bar H$ are randomly generated. 
%Recall that in problem \eqref{relaxation 2}, there are two tuning parameters, $c$ and $\varepsilon$, where $c$ is introduced to tune the tradeoff between the rank constraint and the perturbation cost on $\|K\|_{0,\mathcal{L}}$, and $\varepsilon$ is introduced to guarantee the controllability of the perturbed system.
We first fix $\varepsilon\equiv0.1$ and test the performance with different $c$.
%The results are shown by Table \ref{P0 effect of c} and Fig. \ref{P0 figure}.
From Table \ref{P0 effect of c}, we can see that when $c$ is too small, the perturbed system does not lose a rank.
%After $c$ is larger than a threshold value, the perturbed system begins to lose ranks.
As $c$ increases, the perturbed system has less and less ranks. Row 5 shows that as $c$ increases, the value of $\|\tilde K\|_{0}$ increases. Fig. \ref{P0 figure} shows that after $c$ passing the threshold, ${\rm rank}(\bar D(\tilde z)+F\tilde K)$ has a fast decreasing period (resp. $\|\tilde K\|_{0}$ has a fast increasing period) as $c$ keeps increasing, and after $c$ is larger than another value, ${\rm rank}(\bar D(\tilde z)+F\tilde K)$ decreases (resp. $\|\tilde K\|_{0}$ increases) much slower and tends to constant. Given a perturbed system derived under a specific $c$, we use the mechanism introduced at the second last paragraph of Section \ref{rank relaxation section} to check which data items can be inferred and which cannot. When $c=1$, $(T_1(0),T_7(0),T_8(0),T_9(0),V_{6})$ can be inferred; when $c=1.2$, only $(T_8(0),T_9(0))$ can be inferred; when $c\geq1.5$, no entry can be inferred.
%In this case study, we may take $c\approx2$ as the best choice of $c$, as we desire both non-strong-observability and minimum $\|K\|_{0,\mathcal{L}}$.

%{\color{blue}
%\begin{remark}
%As pointed out in Remark \ref{remark:partial data privacy}, even if a system is not strongly observable, it is possible that the data requester could uniquely determine partial data items of the system. In Section \ref{partial data privacy discussion}, we provide a set of conditions under which specific data items cannot be inferred by the data requester. One can make use of these conditions to check which data items of the perturbed system cannot be inferred by the data requester. For example, in the case of $c=3$ in the above problem, for the perturbed system, only $x_3(0)$ and $x_5(0)$ could be inferred by the data requester, while all the other data items cannot be inferred.\oprocend
%\end{remark}
%}

%A negative message indicated by Fig. \ref{P0 figure} is that the value of $\|K\|_{0,\mathcal{L}}$ around the threshold (around 2) is quite sensitive to the increasing of $c$. To guarantee the non-strong-observability, we may choose $c$ to be a large enough value. However, a marginal

We next verify that the positive semidefinite condition with the introduction of $\varepsilon$ in \eqref{relaxation 2} can guarantee controllability of the perturbed system. We fix $c\equiv2$ and test the cases $\varepsilon=0.01,0.05,0.10,0.50,1.00,5.00$. The perturbed system is controllable for all the tested values.%of $\varepsilon$ in the simulation.

\begin{table}[!t]
\renewcommand{\arraystretch}{1}
\caption{Effect of $c$ with $\varepsilon\equiv0.1$ for problem $\tilde{\mathbb{P}}_{0}$}
\label{P0 effect of c}
\scriptsize
\centering
\begin{tabular}{|l|l|l|l|l|l|}
\hline
$c$ & $0.5$ & $0.8$ & $1.0$ & $2.0$ & $3.0$\\
%\hline
%$\varepsilon$ & $0.1$ & $0.1$ & $0.1$ & $0.1$ & $0.1$\\
%\hline
%S.O. & Yes & No & No & No & No\\
\hline
Controllability & Yes & Yes & Yes & Yes & Yes\\
\hline
${\rm rank}(\bar D(\tilde z)+F\tilde K)$ & $30$ & $27$ & $23$ & $22$ & $17$\\
\hline
$\|\tilde K\|_{0}$ & $28$ & $72$ & $93$ & $300$ & $492$\\
\hline
\end{tabular}
\end{table}
\normalsize

%\begin{table}[!t]
%\renewcommand{\arraystretch}{1.1}
%\caption{Effect of $\varepsilon$ with $c\equiv2$ for problem $\tilde{\mathbb{P}}_{0,\mathcal{L}}$}
%\label{P0 effect of epi}
%\scriptsize
%\centering
%\begin{tabular}{|l|l|l|l|l|l|l|}
%\hline
%%$c$ & $1.0$ & $1.5$ & $2.0$ & $2.5$ & $3.0$\\
%%\hline
%$\varepsilon$ & $0.01$ & $0.05$ & $0.10$ & $0.50$ & $1.00$ & $5.00$\\
%%\hline
%%Strong Observability & No & No & No & No & Yes\\
%\hline
%Controllability & Yes & Yes & Yes & Yes & Yes & Yes\\
%%\hline
%%$n+p-{\rm rank}(D(z)+FK)$ & $19$ & $3$ & $2$ & $1$ & $0$\\
%%\hline
%%$\|K\|_{0,\mathcal{L}}$ & $1234$ & $1236$ & $1245$ & $1249$ & $1255$\\
%\hline
%\end{tabular}
%\end{table}
%\normalsize

%\subsubsection{Problem $\tilde{\mathbb{P}}_{2}$}

%In this case, we also construct the constraint set $\mathbb{K}$ by randomly choosing $100$ positions at which the corresponding entries of $K$ should be $0$. Different from $\mathbb{P}_{0,\mathcal{L}}$, the structure of $\mathbb{K}$ here guarantees that the perturbed system remains controllable under any realization of $K$. We then construct the matrices $\Psi_\kappa$'s and $\Phi_\kappa$'s according to $\mathbb{K}$.

\begin{table}[!t]
\renewcommand{\arraystretch}{1}
\caption{Method of Section \ref{Section P2} for problem $\tilde{\mathbb{P}}_{2}$}
\label{P2 alternative method}
\scriptsize
\centering
\begin{tabular}{|l|l|l|l|l|l|l|}
\hline
%$c$ & $1.0$ & $1.5$ & $2.0$ & $2.5$ & $3.0$\\
%\hline
$\rho$ & $21$ & $16$ & $9$ & $7$ & $5$ & $3$\\
\hline
%Strong Observability & No & No & No & No & Yes\\
%\hline
${\rm rank}(\bar D(\tilde z)+F\tilde K)$ & $20$ & $15$ & $8$ & $6$ & $4$ & $2$\\
\hline
$\|[\bar H,\Pi]\tilde K\|_{2}$ & $0$ & $0.48$ & $1.15$ & $1.23$ & $1.40$ & $1.58$\\
\hline
\end{tabular}
\end{table}
\normalsize

\begin{figure}[H]
\begin{center}
\includegraphics[width=1\linewidth]{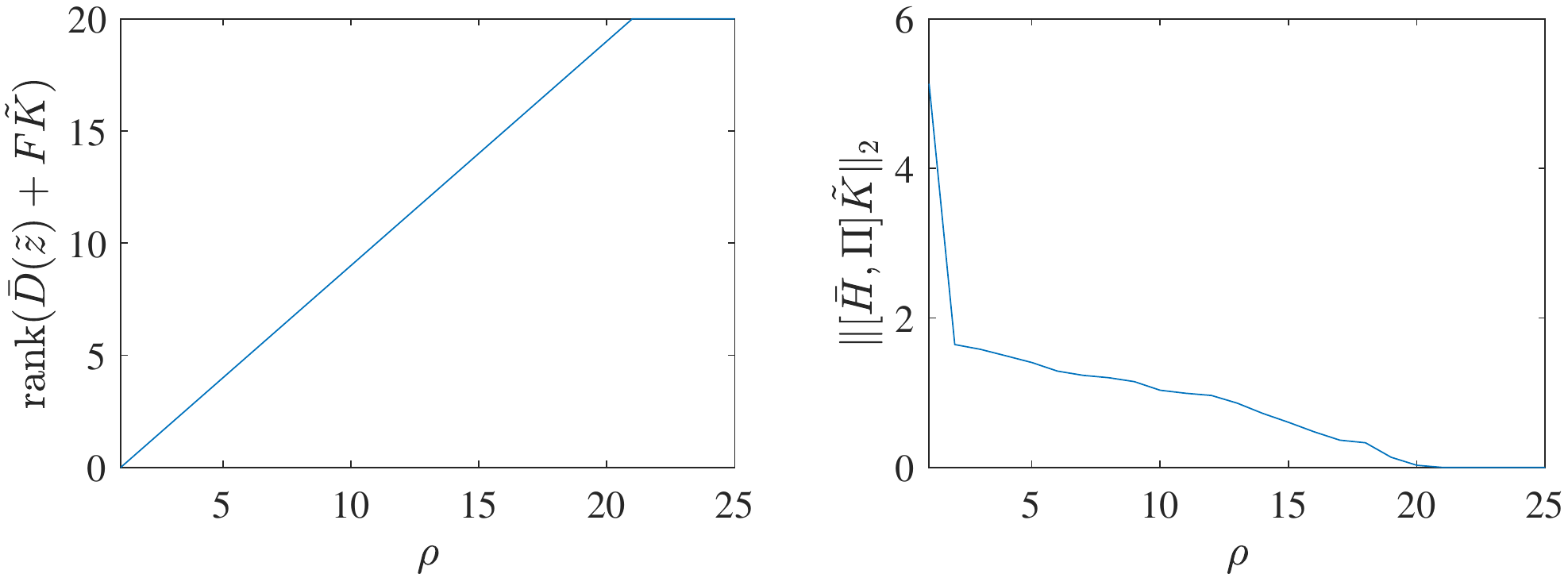}
\caption{Method of Section \ref{Section P2} for problem $\tilde{\mathbb{P}}_{2}$}
\label{P2 figure method of Section V}
\end{center}
\end{figure}

\begin{figure}[h]
\begin{center}
\includegraphics[width=1\linewidth]{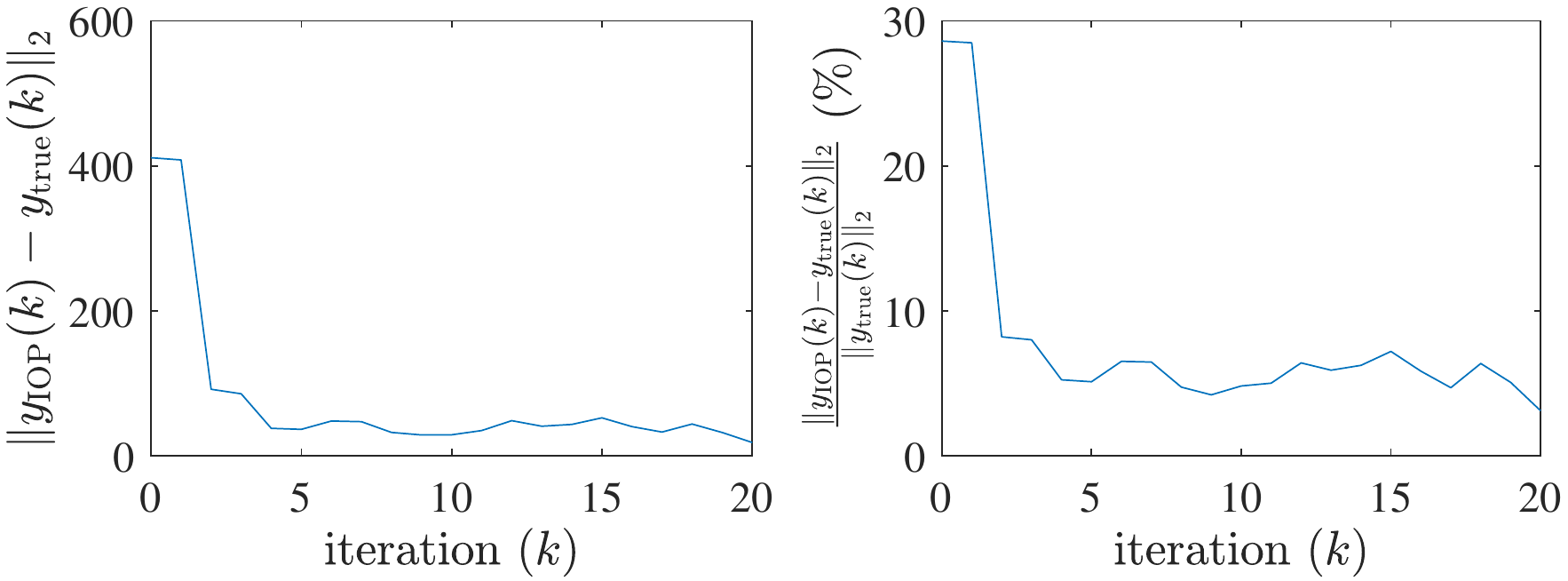}
\caption{Data disutility of problem $\tilde{\mathbb{P}}_{2}$}
\label{P2 data disutility}
\end{center}
\end{figure}

\textbf{Problem $\tilde{\mathbb{P}}_{2}$}. In this case, $\bar G$ and $\bar H^e$ are randomly generated such that Assumption \ref{asm: q equals p} holds.
%We first apply the method of Section \ref{P0 section} given by \eqref{P2 method one}. Notice that for problem $\tilde{\mathbb{P}}_2$ we do not consider the controllability constraint, i.e., the last constraint of \eqref{P2 method one} is dropped. The results are shown by Table \ref{P2 effect of c} and Fig. \ref{P2 figure}. Similar tradeoff as the case of $\tilde{\mathbb{P}}_{0,\mathcal{L}}$ is observed.
We then apply Algorithm \ref{Algorithm P2}. In the simulation, we have $n=10$, $q^e=10$, ${\rm rank}(F)=20$ and ${\rm Tr}(\bar A)=10$. Hence, we have $\tilde z={{\rm Tr}(\bar A)}/{n}=1$.
%Moreover, by Lemma \ref{lemma P2 upper bound}, problem \eqref{fixed z P2} is feasible for all $\rho\geq1$. Given $\rho\geq1$, the perturbation matrix $\tilde K$ is constructed according to Lemma \ref{lemma P2 upper bound} with $z$ taking the value of $\tilde z=1$.
%The simulation results are shown by Table \ref{P2 alternative method} and Fig. \ref{P2 figure method of Section V}.
Table \ref{P2 alternative method} and Fig. \ref{P2 figure method of Section V} show that ${\rm rank}(\bar D(\tilde z)+F\tilde K)=\rho-1<\rho$ for each $\rho$. This verifies that the construction of $K$ given by Lemma \ref{lemma P2 upper bound} is feasible for problem \eqref{fixed z P2}.
%In particular, the result $\|[\bar H,\Pi]\tilde K\|_2=0$ in the range of $\rho\in[21,25]$ verifies (ii) of Lemma \ref{lemma P2 upper bound}.
Table \ref{P2 alternative method} and Fig. \ref{P2 figure method of Section V} also show that the smaller the value of $\rho$, the larger the value of $\|[\bar H,\Pi]\tilde K\|_2$.
%This matches our intuition that a smaller $\rho$ requires to lose more ranks, which is realized by a larger perturbation.
%Comparing Table \ref{P2 effect of c} and Table \ref{P2 alternative method}, we can see that to achieve the same rank of the perturbed system, the $\ell_2$ norm of the perturbation $\tilde K$ generated by the the method of Section \ref{Section P2} is slightly better. %More importantly, to achieve a specific ${\rm rank}(\bar D(\tilde z)+F\tilde K)$, the method of Section \ref{Section P2} does not need empirical tuning of $c$. Instead, one can set $\rho={\rm rank}(\bar D(\tilde z)+F\tilde K)+1$ and analytically construct $\tilde K$ by Lemma \ref{lemma P2 upper bound}.
%Furthermore, for a given $c$, the method of Section \ref{P0 section} needs to solve an SDP, while for a given $\rho$, the method of Section \ref{Section P2} only needs to perform an SVD operation, which is computationally much cheaper. In the simulation, the average time of Table \ref{P2 effect of c} is 8.18 seconds and that of Table \ref{P2 alternative method} is only 0.0048 seconds.
Given a perturbed system derived under a specific $\rho$, we use the mechanism introduced at the second last paragraph of Section \ref{rank relaxation section} to check which data items can be inferred and which cannot. When $\rho=21$, $(T_{10}(0),V_3,V_8)$ can be inferred; when $\rho=19$, only $V_8$ can be inferred; when $\rho\leq18$, for any $i\in\{1,2,\cdots,10\}$, no entry can be inferred. We then design a state feedback controller $u^c$ such that each $x_i$ is stabilized at 21.5 degrees. In the control problem, each $V_i$ is viewed as an external noise and is generated as a random integer between 0 to 10 at each iteration. The data disutility of problem $\tilde{\mathbb{P}}_{2}$ is shown in Fig. \ref{P2 data disutility}, in which $y_{\rm true}$ is the unperturbed output and $y_{\rm IOP}$ is the perturbed output (IOP indicates input-output perturbations). We can see that the data disutility $\|y_{\rm IOP}(k)-y_{\rm true}(k)\|_2$ is below 10\% of $\|y_{\rm true}(k)\|_2$ after 5 iterations.

\begin{figure}[h]
\begin{center}
\includegraphics[width=1\linewidth]{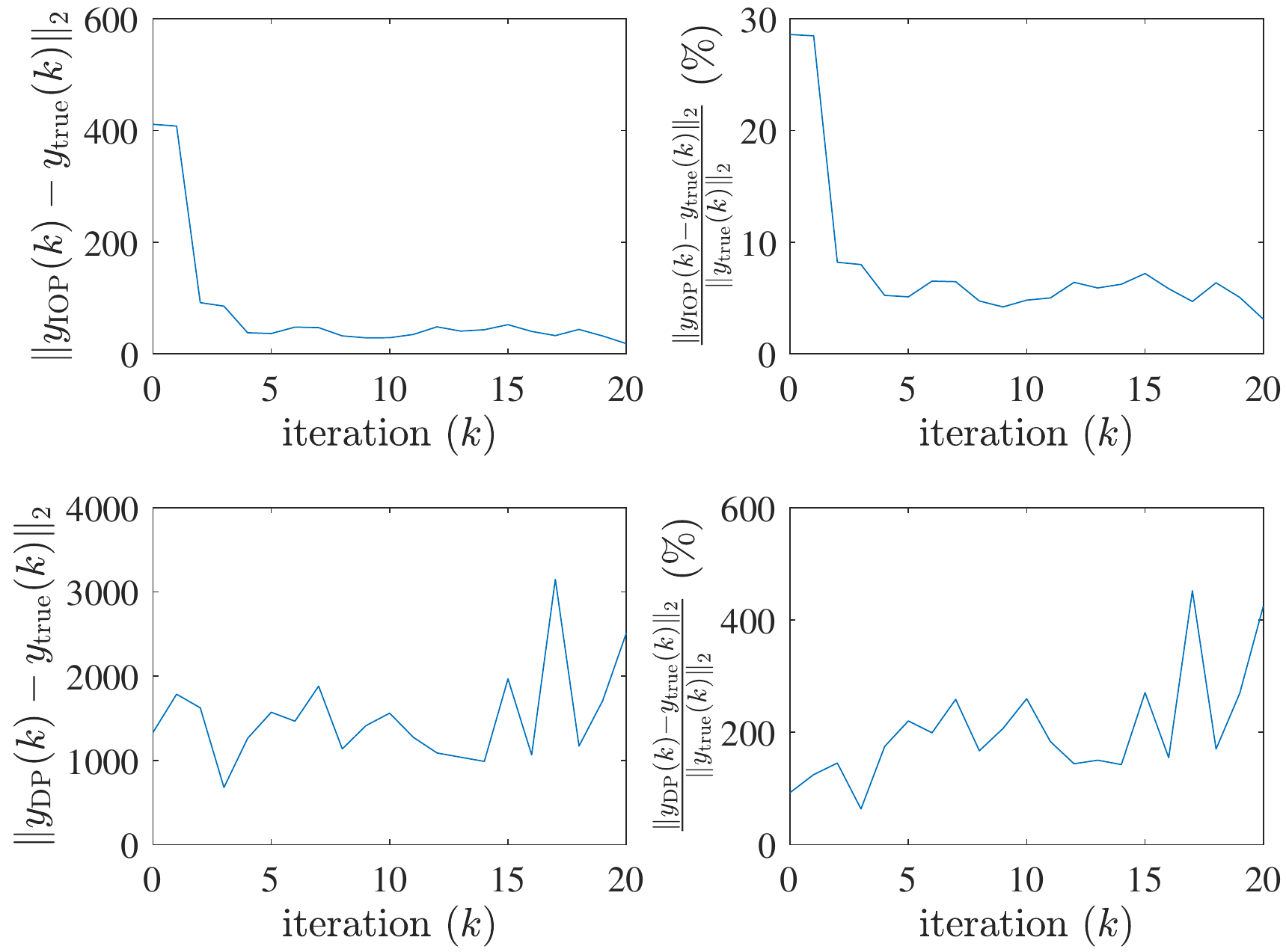}
\caption{Comparison with differential privacy for problem $\tilde{\mathbb{P}}_2$}
\label{P2 DP comparison}
\end{center}
\end{figure}

We also simulate the differentially private scheme in the paper \citep{JLN-GJP:2014} with $\varepsilon=0.1$ (i.e., $0.1$-differential privacy) on the HVAC problem. Fig. \ref{P2 DP comparison} shows the comparison with our algorithm for problem $\tilde{\mathbb{P}}_2$ in terms of data utility. In Fig. \ref{P2 DP comparison}, $y_{\rm true}$ and $y_{\rm IOP}$ have the same meanings as those in Fig. \ref{P2 data disutility}, and $y_{\rm DP}$ is the perturbed output by \citep{JLN-GJP:2014}'s scheme (DP indicates differential privacy). The first row shows data disutility of our method and the second row shows that of the paper \citep{JLN-GJP:2014}. From Fig. \ref{P2 DP comparison}, we can see that our method achieves much better data utility than the differential privacy method of the paper \citep{JLN-GJP:2014} when $\varepsilon=0.1$.
%Moreover, in our method, the perturbations are eventually stable since they are added in a feedback manner. If the steady states are zero, the perturbations will decrease to zero. In contrary, in differential privacy, the perturbations constantly oscillate since they are added in an open-loop manner. This is another advantage of our approach compared with differential privacy. Due to space limitation, Fig. \ref{P2 DP comparison} is not included in the revision.

%We have also simulated the differentially private scheme in the paper \citep{JLN-GJP:2014} with 0.1-differential privacy, and our method achieves much better data utility.

%This verifies that our approach for problem $\mathbb{P}_2$ is reasonable to maintain data utility.}

We also simulate the SDP approach \eqref{P2 method one} with the last constraint dropped. %Similar tradeoff as the case of $\tilde{\mathbb{P}}_{0,\mathcal{L}}$ is observed.
The average time of solving the SDP of \eqref{P2 method one} once is 8.18 seconds, while the average time of running Algorithm \ref{Algorithm P2} once is 0.0048 seconds. As mentioned in the paragraph right below \eqref{P2 method one}, one only needs to tune $\rho$ for at most $\min\{n+q,n+p\}$ times. For this example, $\min\{n+q,n+p\}=20$. Hence, for the worst case, one needs to run Algorithm \ref{Algorithm P2} twenty times and the total running time is approximately $0.0048\times 20=0.096$ seconds, which is still much shorter than solving the SDP of \eqref{P2 method one} once. This verifies that Algorithm \ref{Algorithm P2} is computationally more efficient than the SDP approach.

\section{Conclusions}

This paper formulates the problem of perturbation design to achieve privacy-preserving data release of linear dynamic networks. The computational complexity of the formulated optimization problem is analyzed. An SDP relaxation for the $\ell_0$ minimization is derived. For a class of $\ell_2$ minimizations, we provide a computationally more efficient method which can return an analytic feasible solution. A case study on an HVAC system is conducted to validate the efficacy of the developed techniques.

\section{Appendix\label{introduction to l diversity}}

%\subsection{Introduction to $\ell$-diversity and extension to Definition \ref{def: agent-wise privacy}\label{introduction to l diversity}}
In this section, we first provide an introduction to $\ell$-diversity. After that, we illustrate how to extend $\ell$-diversity to construct Definition \ref{def: agent-wise privacy} in our problem setting.

\begin{table}[h]
\renewcommand{\arraystretch}{1}
\caption{A 3-diverse salary/disease table}
\label{l diversity table}
\scriptsize
\centering
\begin{tabular}{|l|l|l|l|}
\hline
ZIP code & Age & Salary & Disease\\
\hline
476** & 2* & 3K & gastric ulcer\\
476** & 2* & 4K & gastritis\\
476** & 2* & 5K & stomach cancer\\
\hline
4790* & $\geq40$ & 6K & gastritis\\
4790* & $\geq40$ & 11K & flu\\
4790* & $\geq40$ & 8K & brochitis\\
\hline
476** & 3* & 7K & bronchitis\\
476** & 3* & 9K & pneumonia\\
476** & 3* & 10K & stomach cancer\\
\hline
\end{tabular}
\end{table}
\normalsize

Informally speaking, the notion of $\ell$-diversity requires that, given the adversary's observations, there is adequate diversity in each sensitive attribute of the dataset in the released table. The work \cite{NL-TL-SV:2007} formally defines $\ell$-diversity as that each equivalence class of the released table has at least $\ell$ ``well-represented'' values for each sensitive attribute. An equivalence class of an anonymized table is a set of records that share the values of the attributes the adversary may know. Having $\ell$ ``well-represented'' values essentially means that the probabilities of these $\ell$ values are close to each other and meanwhile the total probability of these $\ell$ values is significant, e.g., equal or close to 1. We adopt the following example from \cite{NL-TL-SV:2007} to illustrate the notion of $\ell$-diversity. Table \ref{l diversity table} is an anonymized table with four attributes, namely, ZIP code, Age, Salary and Disease, in which the adversary might observe ZIP codes and Ages of some records, while Salary and Disease are sensitive attributes which should not be disclosed to the adversary. Each * represents an anonymized digit. Each equivalence class shares the values of ZIP code and Age. So there are three equivalence classes: rows 1--3, rows 4--6 and rows 7--9. Since each equivalence class has three different values for each of Salary and Disease, Table \ref{l diversity table} has 3-diversity. Assume that the adversary knows that a specific participant has ZIP code 47630 and age 33 and aims to infer this participant's salary and type of disease. Through Table \ref{l diversity table}, the adversary can tell that this participant's record belongs to the equivalence class formed by the last three rows. However, since this equivalence class has 3-diversity, the adversary cannot uniquely determine the participant's salary or type of disease.

\begin{table}[h]
\renewcommand{\arraystretch}{1}
\caption{An illustrative example}
\label{illustrative table}
\scriptsize
\centering
\begin{tabular}{|l|l|l|}
\hline
$y$ & $x_1(0)$ & $x_2(0)$\\
\hline
2 & 1 & 1\\
2 & 2 & 0\\
2 & 0 & 2\\
$\vdots$ & $\vdots$ & $\vdots$\\
\hline
3 & 1 & 2\\
3 & 2 & 1\\
3 & 3 & 0\\
$\vdots$ & $\vdots$ & $\vdots$\\
\hline
%4 & 1 & 3\\
%4 & 3 & 1\\
%4 & 2 & 2\\
%$\vdots$ & $\vdots$ & $\vdots$\\
%\hline
$\vdots$ & $\vdots$ & $\vdots$\\
\hline
\end{tabular}
\end{table}
\normalsize

\begin{table*}[t]
\renewcommand{\arraystretch}{1}
\caption{Hypothetical released table}
\label{Hypothetical released table}
\scriptsize
\centering
\begin{tabular}{|l|l|l|l|l|l|l|l|l|l|l|l|l|l|}
\hline
$\{y(k)\}$ & $x_1^n(0)$ & $\cdots$ & $x_{d_x^n}^n(0)$ & $\{u_1^n(k)\}$ & $\cdots$ & $\{u_{d_u^n}^n(k)\}$ & $x_1^t(0)$ & $\cdots$ & $x_{d_x^t}^t(0)$ & $\{u_1^t(k)\}$ & $\cdots$ & $\{u_{d_u^t}^t(k)\}$\\
\hline
\multicolumn{13}{|c|}{\multirow{2}{*}{Equivalence class 1} }\\
\multicolumn{13}{ |c|  }{}\\
\hline
\multicolumn{13}{|c|}{\multirow{2}{*}{Equivalence class 2} }\\
\multicolumn{13}{ |c|  }{}\\
\hline
\multicolumn{13}{|c|}{\multirow{2}{*}{$\vdots$} }\\
\multicolumn{13}{ |c|  }{}\\
\hline
\end{tabular}
\end{table*}
\normalsize

The notion of diversity can be interpreted as a measure of uncertainty: a larger diversity indicates that the uncertainty on the sensitive attributes is also larger. In our paper, to make an analogy between $\ell$-diversity and our privacy notion Definition \ref{def: agent-wise privacy}, we can use the output sequence $\{y(k)\}$
%and the system matrices $(\hat A,\hat B,\hat G,\hat H)$ of the perturbed system
as the label of an equivalence class and view each of the adversarial data requester's target entries (i.e., $x^t(0)$ and $u^t$) as a sensitive attribute. To fix idea, we first provide an illustrative example. For simplicity, let $\hat A$, $\hat B$ and $\hat H$ all be zero matrices and $\hat G=[1,1]$. Then the system becomes the single constant output equation $y=x_1(0)+x_2(0)$, where $x_1(0)$ and $x_2(0)$ are sensitive attributes and $y$ is the data requester's observation. The released table is in the form of Table \ref{illustrative table}. In Table \ref{illustrative table}, each real number for $y$ generates an equivalence class and each equivalence class has infinite diversity/uncertainty on both $x_1(0)$ and $x_2(0)$. Using the observed output $y$, say $y=2$, the data requester can determine that $x_1(0)$ and $x_2(0)$ must take values in one record of the equivalence class corresponding to $y=2$. Since $x_1(0)$ and $x_2(0)$ could be any point on the line $x_1(0) + x_2(0) = 2$, the diversity/uncertainty on $x_1(0)$ and $x_2(0)$ is infinite.
%Given its observation $y$ and auxiliary information $(\hat g_1,\hat g_2)$, the data requester can only settle down to a specific equivalence class, say the equivalence with $y=2$ and $(\hat g_1,\hat g_2)=(1,1)$, but has infinite diversity/uncertainty on both $x_1(0)$ and $x_2(0)$. Hence, the privacy is preserved against the data requester's attack.

%\begin{table*}[t]
%\renewcommand{\arraystretch}{1.1}
%\caption{Hypothetical released table}
%\label{Hypothetical released table}
%%\scriptsize
%\centering
%\begin{tabular}{|l|l|l|l|l|l|l|l|l|l|l|l|l|l|}
%\hline
%$\{y(k)\}$ & $x_1^n(0)$ & $\cdots$ & $x_{d_x^n}^n(0)$ & $\{u_1^n(k)\}$ & $\cdots$ & $\{u_{d_u^n}^n(k)\}$ & $x_1^t(0)$ & $\cdots$ & $x_{d_x^t}^t(0)$ & $\{u_1^t(k)\}$ & $\cdots$ & $\{u_{d_u^t}^t(k)\}$\\
%\hline
%\multicolumn{13}{|c|}{\multirow{2}{*}{Equivalence class 1} }\\
%\multicolumn{13}{ |c|  }{}\\
%\hline
%\multicolumn{13}{|c|}{\multirow{2}{*}{Equivalence class 2} }\\
%\multicolumn{13}{ |c|  }{}\\
%\hline
%\multicolumn{13}{|c|}{\multirow{2}{*}{$\vdots$} }\\
%\multicolumn{13}{ |c|  }{}\\
%\hline
%\end{tabular}
%\end{table*}

We next illustrate how to construct the released table for the general case of our problem. The following notations are consistent with those used in Section \ref{privacy notion section}. Our hypothetical released table has the form of Table \ref{Hypothetical released table}.
%\begin{table*}[t]
%\renewcommand{\arraystretch}{1.1}
%\caption{Hypothetical released table}
%\label{Hypothetical released table}
%%\scriptsize
%\centering
%\begin{tabular}{|l|l|l|l|l|l|l|l|l|l|l|l|l|l|}
%\hline
%$\{y(k)\}$ & $x_1^n(0)$ & $\cdots$ & $x_{d_x^n}^n(0)$ & $\{u_1^n(k)\}$ & $\cdots$ & $\{u_{d_u^n}^n(k)\}$ & $x_1^t(0)$ & $\cdots$ & $x_{d_x^t}^t(0)$ & $\{u_1^t(k)\}$ & $\cdots$ & $\{u_{d_u^t}^t(k)\}$\\
%\hline
%\multicolumn{13}{|c|}{\multirow{2}{*}{Equivalence class 1} }\\
%\multicolumn{13}{ |c|  }{}\\
%\hline
%\multicolumn{13}{|c|}{\multirow{2}{*}{Equivalence class 2} }\\
%\multicolumn{13}{ |c|  }{}\\
%\hline
%\multicolumn{13}{|c|}{\multirow{2}{*}{$\vdots$} }\\
%\multicolumn{13}{ |c|  }{}\\
%\hline
%\end{tabular}
%\end{table*}
In Table \ref{Hypothetical released table}, each equivalence class is a set of initial states and inputs which produce the same $\{y(k)\}$. To be more specific, each equivalence class labeled by a specific output sequence $y_{[0,\kappa]}$ includes the target entries $(x^t(0),u^t_{[0,\kappa]})$ of $\Delta_{\hat A,\hat B,\hat G,\hat H}(y_{[0,\kappa]})$, together with admissible non-target entries $(x^n(0),u^n_{[0,\kappa]})$.
%each feasible output sequence $\{y(k)\}$ of system $(\hat A,\hat B,\hat G,\hat H)$ specifies an equivalence class such that each element includes particular initial state and input sequence which produce output $\{y(k)\}$, that is, $y(k)=\hat G\hat A^kx(0)+\sum_{m=0}^{k-1}\hat G\hat A^{k-1-m}\hat Bu(m)+\hat Hu(k)$ for any $k\in\mathbb{N}$, with $x(0)$ being the composition of $x^t(0)$ and $x^n(0)$, and $u(k)$ being the composition of $u^t(k)$ and $u^n(k)$ for any $k\in\mathbb{N}$.
Similar to $\ell$-diversity, our privacy goal is to guarantee that each equivalence class has adequate diversity/uncertainty on each of the target entries $(x_1^t(0),\cdots,x_{d_x^t}^t(0),u_1^t,\cdots,u_{d_u^t}^t)$. In $\ell$-diversity, sensitive attributes take discrete values and the diversity on each sensitive attribute is defined by the number of different valuations for that attribute in each equivalence class. In contrast, the target entries $x^t(0)$ and $u^t$ in our paper are continuous-valued and thus cannot be enumerated (i.e., uncountable). Hence, we need to introduce a new measure to quantify the diversity/uncertainty. In this paper, we propose to measure the diversity/uncertainty by the diameter of the set $\Delta_{\hat A,\hat B,\hat G,\hat H}(y_{[0,\kappa]})$. For each target entry, a larger diameter indicates a larger range of admissible valuations and thus a larger diversity/uncertainty. An infinite diameter achieves the largest possible diversity/uncertainty. Hence, for our problem setting, we say that privacy is preserved if the diameter of the set $\Delta_{\hat A,\hat B,\hat G,\hat H}(y_{[0,\kappa]})$ is infinite for any feasible output sequence $y_{[0,\kappa]}$ for any $\kappa\in\mathbb{N}$.
%we consider the case of infinite diversity/uncertainty to be adequate, for which we say that privacy is preserved;
Please refer to Section \ref{privacy notion section} for the detailed definitions and discussions.

\bibliographystyle{agsm-nq}
\bibliography{MZ,YL}

\end{document}